\documentstyle[amstex,amssymb]{amsart}
\numberwithin{equation}{section}
\newtheorem{thm}{Theorem}[section] 
\newtheorem{cor}[thm]{Corollary}
\newtheorem{prop}[thm]{Proposition}
\newtheorem{lem}[thm]{Lemma}
\newcommand{\bv}{\mbox{\bf v}}
\newcommand{\comment}[1]{}
\address{
Department of Mathematics, Kobe University, 
Rokko, Kobe 657, Japan
}
\email{
noumi@@math.s.kobe-u.ac.jp\quad yamaday@@math.s.kobe-u.ac.jp
}
\catcode`\@11 \newdimen\b@xsize
\newtoks\p@rtition
\newcount\r@
\newtoks\@nds
\newif\if@cont
\def\Young#1(#2){\def\aho{#1}\ifx\aho\empty\b@xsize6pt\else\b@xsize\aho\fi
\global\p@rtition{#2,\aho}\global\@conttrue\expandafter
\getnum\the\p@rtition
\ifnum\r@>0\global\setbox1\hbox to\b@xsize{\hfil
\vrule width0ptheight\b@xsize depth0pt}
\vcenter{\hbox{\vbox{\offinterlineskip\tabskip0pt\def\@nd{&}
\halign{\vrule height\b@xsize##&&\copy1##&\vrule##\cr\@@hline\thisline@\cr}}}}
\else\emptyset\fi}%
\def\young#1(#2){\def\aho{#1}\ifx\aho\empty\b@xsize2.5pt\else\b@xsize\aho\fi
\global\p@rtition{#2,\aho}\global\@conttrue\expandafter
\getnum\the\p@rtition
\ifnum\r@>0\global\setbox1\hbox to\b@xsize{\hfil
\vrule width0ptheight\b@xsize depth0pt}\hbox{\vtop{\offinterlineskip\tabskip0pt
\def\@nd{&}\halign{\vrule height\b@xsize##&&\copy1##&\vrule##\cr
\@@hline\thisline@\cr}}}\else\emptyset\fi}%
\def\thisline@{\cr\r@wofbox\cr\@@hline\if@cont\let\noriP
\tr@ns\else\let\noriP\relax\fi\noriP}%
\def\tr@ns{\expandafter\getnum\the\p@rtition\thisline@}%
\def\getnum#1,#2\aho{\def\aho{#2}\ifx\aho\empty\def\aho{0}\global\@contfalse\fi
\ifnum#1=0\global\@contfalse\fi
\global\r@#1\global\multiply\r@\tw@
\global\p@rtition{#2\aho}}%
\def\@@hline{\omit\mscount\r@\loop\ifnum\mscount>\@ne\sp@n\repeat\hrulefill&}
\def\r@wofbox{\global\@nds{}\global\count\@ne\r@\gdef\foolish{\global\@nds
\expandafter{\the\@nds\@nd}\global\advance\count\@ne -\@ne
\ifnum\count\@ne>0}\boke\the\@nds}%
\def\boke{\foolish\global\let\ahoaho\boke\else\global\let\ahoaho\relax\fi\ahoaho}%
\catcode`\@\active
\def\FigOka{
\begin{figure}
\setlength{\unitlength}{1.1mm}
\begin{picture}(100,85)(-10,-10)
\hbox{
\hskip4mm
\put(0,0){\young(2)}
\put(20,0){\young(1)}
\put(40,0){\young(1,1)}
\put(60,0){\young(3,1,1)}
\put(80,0){\young(5,3,1,1)}
\put(10,15){$\phi$} \put(13,17){$Q_0$}\put(13,12){\small$(0,0)$}
\put(30,15){$\phi$} \put(33,17){$Q_1$}\put(33,12){\small$(1,0)$}
\put(50,15){\young(2)} \put(53,17){$Q_2$}
\put(70,15){\young(4,2)} \put(73,17){$Q_3$}
\put(90,15){\young(6,4,2)} \put(93,17){$Q_4$}
\put(0,30){\young(1)}\put(3,32){$R_{-1}$}\put(3,27){\small$(0,1)$}
\put(20,30){$\phi$} \put(23,32){$R_0$}\put(23,27){\small$(1,1)$}
\put(40,30){\young(1)} \put(43,32){$R_1$}
\put(60,30){\young(3,1)} \put(63,32){$R_2$}
\put(80,30){\young(5,3,1)} \put(83,32){$R_3$}
\put(10,45){\young(2)}
\put(30,45){\young(1,1)}
\put(50,45){\young(2,1,1)}
\put(70,45){\young(4,2,1,1)}
\put(90,45){\young(6,4,2,1,1)}
\put(0,60){\young(4,2)}
\put(20,60){\young(3,1,1)}
\put(40,60){\young(2,2,1,1)}
\put(60,60){\young(3,2,2,1,1)}
\put(80,60){\young(5,3,2,2,1,1)}
\put(10,75){\young(5,3,1,1)}
\put(30,75){\young(4,2,2,1,1)}
\put(50,75){\young(3,3,2,2,1,1)}
\put(70,75){\young(4,3,3,2,2,1,1)}
\put(90,75){\young(6,4,3,3,2,2,1,1)}
}
\multiput(5,0)(20,0){4}{\vector(1,0){10}}
\multiput(15,15)(20,0){4}{\vector(1,0){10}}
\multiput(5,30)(20,0){4}{\vector(1,0){10}}
\multiput(15,45)(20,0){4}{\vector(1,0){10}}
\multiput(5,60)(20,0){4}{\vector(1,0){10}}
\multiput(15,75)(20,0){4}{\vector(1,0){10}}
\multiput(18,-2)(20,0){4}{\vector(-2,-3){6}}
\multiput(8,13)(20,0){5}{\vector(-2,-3){6}}
\multiput(18,28)(20,0){4}{\vector(-2,-3){6}}
\multiput(8,43)(20,0){5}{\vector(-2,-3){6}}
\multiput(18,58)(20,0){4}{\vector(-2,-3){6}}
\multiput(8,73)(20,0){5}{\vector(-2,-3){6}}
\multiput(8,-12)(20,0){5}{\vector(-2,3){6}}
\multiput(18, 3)(20,0){4}{\vector(-2,3){6}}
\multiput( 8,18)(20,0){5}{\vector(-2,3){6}}
\multiput(18,33)(20,0){4}{\vector(-2,3){6}}
\multiput( 8,48)(20,0){5}{\vector(-2,3){6}}
\multiput(18,63)(20,0){4}{\vector(-2,3){6}}
\end{picture}
\caption{Okamoto polynomials on the $A_2$-lattice}
\label{fig:okamoto}
\end{figure}
}
\def\ExOka{
\par\noindent
{\bf Okamoto polynomials.}
\quad 
In the following, we use the notation 
$Q_\lambda(x)=Q_{m,n}(x)$ for the Okamoto polynomial 
associated with the 3-reduced partition $\lambda=\lambda(m,n)$. 
We give below some examples of Okamoto polynomials 
$Q_\lambda(x)$. 
{\allowdisplaybreaks
\begin{align}
&Q_{(0)}=Q_1=1, \quad
Q_{(1)}=R_1=x, \quad
Q_{(2)}=Q_2=1 + x^2, \notag\\
&Q_{(1, 1)}=-1 + x^2, \quad
Q_{(3, 1)}=R_2=-1 + 2 x^2 + x^4, \notag\\
&Q_{(2, 1, 1)}=-1 - 2 x^2 + x^4, \quad
Q_{(3, 1, 1)}=-5 x + x^5,\notag\\
&Q_{(4, 2)}=Q_3=5 + 5 x^2 + 5 x^4 + x^6, \quad
Q_{(2, 2, 1, 1)}=-5 + 5 x^2 - 5 x^4 + x^6, \notag \\
&Q_{(4, 2, 1, 1)}=-7 - 14 x^4 + x^8, \quad
Q_{(5, 3, 1)}=R_3=-35 x + 14 x^5 + 8 x^7 + x^9, \notag \\
&Q_{(5, 3, 1, 1)}=
  25 - 75 x^2 - 50 x^4 - 10 x^6 + 5 x^8 + x^{10}, \notag \\
&Q_{(6, 4, 2)}=Q_4=
  175 + 350 x^2 + 175 x^4 + 140 x^6 + 65 x^8 + 
   14 x^{10} + x^{12}, \notag \\
&Q_{(7, 5, 3, 1)}=R_4=
  1225 - 4900 x^2 - 4900 x^4 - 980 x^6 \notag\\
&\quad + 350 x^8 + 420 x^{10} + 140 x^{12} 
+ 20 x^{14} + x^{16}, \notag \\
&Q_{(8, 6, 4, 2)}=Q_5=
  67375 + 134750 x^2 + 202125 x^4 + 107800 x^6  \notag \\
 &\quad+  42350 x^8 + 20020 x^{10} + 8050 x^{12} + 2200 x^{14} + 
   355 x^{16} + 30 x^{18} + x^{20}. \notag
\end{align}
}
Note that the original Okamoto polynomials are given by 
\begin{equation}
\begin{align}
&Q_n=Q_{(2n-2,2n-4,\ldots,4,2)}\quad
(n>0), \quad 
Q_{-n}=Q_{(n,n,\ldots,2,2,1,1)}\quad(n\ge 0),\notag\\
&R_n=Q_{(2n-1,2n-3,\ldots,3,1)}\quad( n>0),\quad 
R_{-n}=Q_{(n,n+1,n+1,\ldots,1,1)}\quad(n\ge 0).\notag
\end{align}
\end{equation}
}
\def\FigHer{
\begin{figure}
\setlength{\unitlength}{1.1mm}
\begin{picture}(100,70)(-5,-5)
\hbox{
\hskip4mm
\put(0,0){$\phi$} \put(3,2){$H_{0,0}$}
\put(20,0){$\phi$} \put(23,2){$H_{1,0}$}
\put(40,0){$\phi$}
\put(60,0){$\phi$}
\put(80,0){$\phi$}
\put(10,15){$\phi$} \put(13,17){$H_{0,1}$}
\put(30,15){\young(1)} \put(33,17){$H_{1,1}$}
\put(50,15){\young(1,1)}\put(53,17){$H_{2,1}$}
\put(70,15){\young(1,1,1)}\put(73,17){$H_{3,1}$}
\put(90,15){\young(1,1,1,1)}\put(93,17){$H_{4,1}$}
\put(20,30){$\phi$}
\put(40,30){\young(2)} \put(44,33){$H_{1,2}$}
\put(60,30){\young(2,2)}
\put(80,30){\young(2,2,2)}
\put(30,45){$\phi$} 
\put(50,45){\young(3)} \put(54,48){$H_{1,3}$}
\put(70,45){\young(3,3)} 
\put(90,45){\young(3,3,3)} 
\put(40,60){$\phi$}
\put(60,60){\young(4)} \put(64,63){$H_{1,4}$}
\put(80,60){\young(4,4)}
}
\multiput(5,0)(20,0){5}{\vector(1,0){10}}
\multiput(15,15)(20,0){4}{\vector(1,0){10}}
\multiput(25,30)(20,0){4}{\vector(1,0){10}}
\multiput(35,45)(20,0){3}{\vector(1,0){10}}
\multiput(45,60)(20,0){3}{\vector(1,0){10}}
\multiput(8,13)(20,0){5}{\vector(-2,-3){6}}
\multiput(18,28)(20,0){4}{\vector(-2,-3){6}}
\multiput(28,43)(20,0){4}{\vector(-2,-3){6}}
\multiput(38,58)(20,0){3}{\vector(-2,-3){6}}
\multiput(18, 3)(20,0){4}{\vector(-2,3){6}}
\multiput(28,18)(20,0){4}{\vector(-2,3){6}}
\multiput(38,33)(20,0){3}{\vector(-2,3){6}}
\multiput(48,48)(20,0){3}{\vector(-2,3){6}}
\end{picture}
\caption{
Generalized Hermite polynomials on the $A_2$-lattice 
}
\label{fig:hermite}
\end{figure}
}
\def\ExHer{
\par\medskip\noindent
{\bf Generalized Hermite polynomials.} \quad
The polynomials $H_{n,1}(x)$ and 
$H_{1,n}(x)$ coincide with the ordinary Hermite polynomials up to 
rescaling.  
{\allowdisplaybreaks
\begin{align}
&H_{0,0}=1, \quad
H_{1,0}=1, \quad
H_{2,0}=3, \quad
H_{3,0}=2^1 3^3, \quad
H_{4,0}=2^2 3^7, \notag \\ &
H_{0,1}=1, \quad
H_{1,1}=3 x, \quad
H_{2,1}=3^3 (-{1 \over 3}+x^2), \notag \\&
H_{3,1}=2^1 3^6(-x+x^3), \quad
H_{4,1}=2^2 3^{11} ({1 \over 3}-2 x^2+x^4), \notag \\ &
H_{0,2}=-3, \quad
H_{1,2}=-3^3({1 \over 3}+x^2),\quad
H_{2,2}=-3^6 ({1 \over 3}+x^4), \notag \\ &
H_{3,2}=-2^1 3^{10}({1 \over 3}+x^2-x^4+x^6), \quad
H_{4,2}=-2^2 3^{16}({5 \over 9}+{10 \over 3}x^4-
{8 \over 3}x^6+x^8), \notag \\ &
H_{0,3}=-2^1 3^3, \quad
H_{1,3}=-2^1 3^6(x+x^3), \notag \\ &
H_{2,3}=-2^1 3^{10}(-{1 \over 3}+x^2+x^4+x^6), \quad
H_{3,3}=-2^2 3^{15}({-5 \over 3}x+2 x^5+x^9), \notag \\ &
H_{4,3}=-2^3 3^{22} ({25 \over 27}-{50 \over 9}x^2-{25 \over 9}x^4-
{20 \over 9}x^6+5 x^8-2 x^{10}+x^{12}), \notag 
\\ &
H_{0,4}=2^2 3^7, \quad
H_{1,4}=2^2 3^{11}({1 \over 3}+2x^2+x^4), \notag \\ &
H_{2,4}=2^2 3^{16}({5\over9}+{10\over3}x^4+{8\over3}x^6+x^8), \notag \\ &
H_{3,4}=2^3 3^{22}({25 \over 27}+{50 \over 9}x^2-{25 \over 9}x^4+
{20 \over 9}x^6+5 x^8+2 x^{10}+x^{12}), \notag \\ &
H_{4,4}=2^4 3^{30} ({875 \over 243}+{3500 \over 81} x^4-{50 \over 9}x^8+
{20 \over 3}x^{12}+x^{16}). \notag
\end{align}
}
}
\def\Lat{
\begin{figure}
\setlength{\unitlength}{1.1mm}
\begin{picture}(100,80)(-15,-15)
\hbox{
\hskip5mm
\put(-1,-1){$\bullet$}
\put(19,-1){$\bullet$}
\put(39,-1){$\bullet$}
\put(59,-1){$\bullet$}
\put(79,-1){$\bullet$}
\put(9,14){$\bullet$}
\put(29,14){$\bullet$}
\put(49,14){$\bullet$}\put(49,18){\small $F$}
\put(69,14){$\bullet$}
\put(89,14){$\bullet$}
\put(-1,29){$\bullet$}
\put(19,29){$\bullet$}
\put(39,29){$\bullet$}\put(42,31){\small $A$}
\put(59,29){$\bullet$}\put(55,31){\small $B$}
\put(79,29){$\bullet$}
\put(9,44){$\bullet$}
\put(29,44){$\bullet$}\put(33,42){\small $E$}
\put(49,44){$\bullet$}\put(49,41){\small $C$}
\put(69,44){$\bullet$}\put(64,42){\small $D$}
\put(89,44){$\bullet$}
\put(21.5,1){\small $O$}
\put(29,-2){\small $\mbox{\bf e}_1$}
\put(12,6){\small $\mbox{\bf e}_2$}
\put(29,4){${\cal C}$}
\put(35,1.5){\small $\Lambda_1$}
\put(28.5,10.5){\small $\Lambda_2$}
}
\multiput(5,0)(20,0){4}{\vector(1,0){10}}
\multiput(15,15)(20,0){4}{\vector(1,0){10}}
\multiput(5,30)(20,0){4}{\vector(1,0){10}}
\multiput(15,45)(20,0){4}{\vector(1,0){10}}
\multiput(18,-2)(20,0){4}{\vector(-2,-3){6}}
\multiput(8,13)(20,0){5}{\vector(-2,-3){6}}
\multiput(18,28)(20,0){4}{\vector(-2,-3){6}}
\multiput(8,43)(20,0){5}{\vector(-2,-3){6}}
\multiput(18,58)(20,0){4}{\vector(-2,-3){6}}
\multiput(8,-12)(20,0){5}{\vector(-2,3){6}}
\multiput(18, 3)(20,0){4}{\vector(-2,3){6}}
\multiput( 8,18)(20,0){5}{\vector(-2,3){6}}
\multiput(18,33)(20,0){4}{\vector(-2,3){6}}
\multiput( 8,48)(20,0){5}{\vector(-2,3){6}}
\end{picture}
\caption{The $A_2$-lattice}
\label{fig:lattice0}
\end{figure}
}
\begin{document}
%%%%%%%%%%%%%%%%%%%%%%%%%%%%%%%%%%%%%%%%%
% Title, authors,...
%%%%%%%%%%%%%%%%%%%%%%%%%%%%%%%%%%%%%%%%%
\par\bigskip
\begin{center}
{\Large \bf
Symmetries in the fourth  Painlev\'e  equation \\ 
and Okamoto polynomials
}
\end{center}
\par\medskip
\begin{center}
Masatoshi NOUMI and Yasuhiko YAMADA \\
{\small
Department of Mathematics, Kobe University
}
\end{center}
\par\bigskip
%%%%%%%%%%%%%%%%%%%%%%%%%%%%%%%%%%%%%%%%%
%\section*{Introduction}
%%%%%%%%%%%%%%%%%%%%%%%%%%%%%%%%%%%%%%%%%
\par\bigskip
It is known by K.\,Okamoto \cite{Ok} that the fourth Painlev\'e equation has
symmetries  under the affine Weyl group of type $A^{(1)}_2$.  
In this paper we propose a new representation of the fourth Painlev\'e 
equation in which the $A^{(1)}_2$-symmetries become clearly visible. 
By means of this representation, we clarify the internal relation between 
the fourth Painlev\'e equation  and the modified KP hierarchy.  
We obtain in particular a complete description of the rational 
solutions of the fourth  Painlev\'e equation in terms of Schur functions. 
This implies that the so-called {\em Okamoto polynomials}, 
which arise from the $\tau$-functions for rational solutions, 
are in fact expressible by the 3-reduced Schur functions.
\footnote{
After completing this paper, the authors were informed by 
K.\,Kajiwara and Y.\,Ohta that they obtained independently 
the expression of Okamato polynomials in terms of Schur functions. }
\par\medskip
%%%%%%%%%%%%%%%%%%%%%%%%%%%%%%%%%%%%%%%%%%%%
\section{A symmetric form of the fourth Painlev\'e equation}
%%%%%%%%%%%%%%%%%%%%%%%%%%%%%%%%%%%%%%%%%%%%
The fourth Painlev\'e equation $\mbox{P}_{\mbox{\small IV}}$ is the following
second order ordinary differential equation 
\begin{equation}
y''={1 \over 2 y}(y')^2+{3 \over 2} y^3+4 t y^2+2(t^2-a) y+
{b \over y}
\label{p4y}
\end{equation}
for the unknown function $y=y(t)$, 
where $'=d/dt$ and $a, b \in {\Bbb C}$ are 
parameters. 
It is known by K.\,Okamoto \cite{Ok} that equation \eqref{p4y} is represented 
as the following  system for the two unknown functions $q=y$ and $p$: 
\begin{equation}
\begin{align}
\label{p4h}
&q'=q(2p-q-2t)-2(v_1-v_2), \\
&p'=p(2q-p+2t)+2(v_2-v_3). \notag 
\end{align}
\end{equation}
This equation, called $\mbox{H}_{\mbox{\small IV}}$, 
is in fact a Hamiltonian system
\begin{equation}
q'={\partial H \over \partial p}, \quad
p'=-{\partial H \over \partial q}
\end{equation}
with polynomial Hamiltonian
\begin{equation}
H= q p^2-q^2 p-2 tpq-2(v_1-v_2)p-2(v_2-v_3) q. 
\label{hamiltonian}
\end{equation}
The parameters $\bv=(v_1,v_2,v_3)$, ($v_1+v_2+v_3=0$) in \eqref{p4h} and $(a,
b)$ in \eqref{p4y} are related through the formulas 
\begin{equation}
%a=1-v_1-v_2+2 v_3,
a=1+3 v_3,\quad b=-2(v_1-v_2)^2. 
\label{parameter-ab}
\end{equation}
The equivalence between \eqref{p4y} and \eqref{p4h} can be checked 
directly, but it requires a tedious calculation. 
(This calculation is fairly simplified by the ``symmetric'' representation 
which we will propose in this paper. 
See the proof of Theorem \ref{symmform} below.)
\par
It is clearly seen from \eqref{p4h} that, if $v_1-v_2=0$ or $v_2-v_3=0$, 
the Hamiltonian system $\mbox{H}_{\mbox{\small IV}}$ 
has {\em classical solutions} such that 
$q=0$ or $p=0$. 
In these cases,  equation \eqref{p4h} is reduced to the Riccati equations 
$p'=-p^2+2tp+2(v_2-v_3)$ and $q'=-q^2-2t-2(v_1-v_2)$ respectively, and 
they are furthermore linearized to Hermite-Weber equations. 
In this sense,  the Hamiltonian system 
$\mbox{H}_{\mbox{\small IV}}$
\eqref{p4h} has  {\em invariant divisors} $q=0$ and $p=0$ along the lines 
$v_1-v_2=0$ and $v_2-v_3=0$, respectively. 
It should be noted that equation \eqref{p4h} has one more typical invariant 
divisor $q-p+2t=0$ along the line $v_1-v_3=1$.  
In fact equation \eqref{p4h} implies
\begin{equation}
\label{iv0}
(q-p+2t)'=-(q-p-2t)(q+p) +2(1-v_1+v_3).  
\end{equation}
It is known by \cite{NO} that these three polynomials $q$, $p$ and $q-p+2t$ 
%together with their B\"acklund transformations, 
generate essentially all the invariant divisors of the fourth Painlev\'e
equation \eqref{p4h}.  
Note that the three simple affine roots 
$1-v_1+v_3$, $v_1-v_2$, $v_2-v_3$ 
of type $A^{(1)}_2$ are already involved in these equations.  
We denote by 
\begin{equation}
V=\{\bv=(v_1,v_2,v_3)\in {\Bbb C}^3\ ; \  v_1+v_2+v_3=0 \}
\end{equation}
the parameter space for the system \eqref{p4h}.
\par\medskip
We now propose to treat the three typical invariant divisors 
$q$, $p$ and $q-p+2t$ 
equally so as to obtain a ``symmetric'' representation of the 
fourth Painlev\'e equation. 
We introduce the three dependent variables 
$f=(f_0,f_1,f_2)$ as follows. 
Fixing a nonzero complex number $c\in{\Bbb C}^\times$, set
\begin{equation} 
\label{deffa}
f_0=c(q-p+2t),\quad f_1=-c q, \quad f_2=c p, 
\end{equation} 
and rescale the independent variable as $x=-t/c$. 
Then we have
\begin{equation}\label{p4c}
\begin{align}
&f_0'=f_0(f_2-f_1)-2c^2(1-v_1+v_3),\notag\\
&f_1'=f_1(f_0-f_2)-2c^2(v_1-v_2),\\
&f_2'=f_2(f_1-f_0)-2c^2(v_2-v_3),\notag
\end{align}
\end{equation}
where $'=d/dx$.  
With the normalization $c=\sqrt{-3/2}$, we set
\begin{equation}
\alpha_0=3(1-v_1+v_3),\quad \alpha_1=3(v_1-v_2),\quad \alpha_2=3(v_2-v_3). 
\label{parameter-alpha}
\end{equation}
Then we have 
\begin{thm}\label{symmform}
The fourth Painlev\'e equation \eqref{p4y} $($ or \eqref{p4h} $)$  
can be written in the following symmetric form: 
\begin{equation}
\begin{align}
\label{p4f}
&f'_0+f_0(f_1-f_2)=\alpha_0, \notag\\
&f'_1+f_1(f_2-f_0)=\alpha_1,  \\
&f'_2+f_2(f_0-f_1)=\alpha_2, \notag
\end{align}
\end{equation}
with normalization
\begin{equation}\label{p4fnorm}
f_0+f_1+f_2=3x, 
\end{equation}
where $'=d/dx$ and $\alpha_0, \alpha_1, \alpha_2\in{\Bbb C}$ are parameters
with $\alpha_0+\alpha_1+\alpha_2=3$. 
\end{thm}
\medskip\noindent
{\em Proof.}\quad
The equation \eqref{p4f} has been derived from the Hamiltonian 
system \eqref{p4h}; it is clear that these two are equivalent.
We will show the equivalence of \eqref{p4y}
and \eqref{p4f} (with normalization \eqref{p4fnorm}).
This gives in fact an easier way to establish the equivalence
between \eqref{p4y} and \eqref{p4h}.
Taking a derivative of the second equation of \eqref{p4f}, we have
\begin{equation}
f_1''+f_1'(f_2-f_0)+f_1(f_2'-f_0')=0.
\notag
\end{equation}
Substituting the first and the third equations of \eqref{p4f} to this, 
we obtain
\begin{equation}
f_1''+f_1'(f_2-f_0)-2f_0f_2f_1+(\alpha_0-\alpha_2)f_1+(f_2+f_0)f_1^2=0.
\notag
\end{equation}
Then, by using the relations
\begin{equation}
f_2-f_0={\alpha_1-f_1' \over f_1}, \quad
f_2+f_0=3x-f_1,
\quad
4f_0f_2=(f_2+f_0)^2-(f_2-f_0)^2,
\notag
\end{equation}
we have
\begin{equation}
f_1''-{1 \over 2} {{f_1'}^2 \over f_1}-{3 \over 2} f_1^3+6x f_1^2+
\left( -{9 \over 2}x^2+(\alpha_0-\alpha_2) \right) f_1+{\alpha_1^2 \over 2}
{1 \over f_1}=0.
\notag
\end{equation}
This is transformed into the equation \eqref{p4y} by the rescaling
$f_1=-cy$, $x=-t/c$, $c=\sqrt{-3/2}$ and the change of parameters 
\eqref{parameter-ab}, \eqref{parameter-alpha}.
$\qed$

\par\medskip
We remark that our equation \eqref{p4f}  
has the following rational solutions:
\begin{equation}
\begin{align} \label{ratsol}
\mbox{(A)} &\qquad(\alpha_0,\alpha_1,\alpha_2; f_0,f_1,f_2) 
= (1,1,1;\, x,x,x),  \\
\mbox{(B)} &\qquad(\alpha_0,\alpha_1,\alpha_2; f_0,f_1,f_2) 
= (3,0,0; 3x,0,0).\notag
\end{align}
\end{equation}
{}From the work of Y.\,Murata \cite{M}, it follows 
that all the rational solutions of \eqref{p4f} are obtained from 
these two particular solutions by B\"acklund transformations. 
% ; other solutions are known to be non-algebraic.
There are classical solutions obtained as B\"acklund transformations 
from the solutions of Riccati type along the three lines $\alpha_0=0$, 
$\alpha_1=0$, $\alpha_2=0$. 
Any other solutions are {\em non-classical} in the sense of H. Umemura 
\cite{U} (see also \cite{NO}, \cite{Ok}).
It should also be noted that our equation \eqref{p4c} reduces to
the Kac-Moerbeke integrable system \cite{KM} in the degenerate limit
$c \rightarrow 0$.

\par\medskip
%%%%%%%%%%%%%%%%%%%%%%%%%%%%%%%%%%%%%%%%%%%%
\section{B\"acklund transformations and the affine Weyl group}
%%%%%%%%%%%%%%%%%%%%%%%%%%%%%%%%%%%%%%%%%%%%

We now discuss symmetries in  the fourth  Painlev\'e equation represented 
by \eqref{p4f} with the normalization of \eqref{p4fnorm}. 
In what follows, we regard $\alpha_0$, $\alpha_1$, $\alpha_2$ 
as coordinate functions (with  $\alpha_0+\alpha_1+\alpha_2=3$ ) 
of the parameter space $V$. 

We consider the affine Weyl group 
$W=\langle s_0, s_1,s_2\rangle$
of type $A^{(1)}_2$ with fundamental relations
\begin{equation}
s_i^2=1,\quad s_i s_{i+1} s_i = s_{i+1} s_i s_{i+1}\quad (i=0,1,2). 
\end{equation}
Here the subscripts are understood as elements of ${\Bbb Z}/3{\Bbb Z}$.  
This convention for subscripts will be applied 
to other variables $\alpha_i$, $f_i$, etc., as well. 
We denote by $\widetilde{W}=\langle s_0,s_1,s_2,\pi \rangle$ 
the extension of $W$  obtained by adjoining the following Dynkin diagram 
automorphism $\pi$:  
\begin{equation}
\pi^3=1,\quad \pi s_i=s_{i+1}\pi \quad(i=0,1,2). 
\end{equation} 
The affine Weyl group $\widetilde{W}$ acts naturally on the coordinate ring 
${\Bbb C}[\alpha]$ of $V$ through the algebra automorphism 
$s_0,s_1,s_2$ and $\pi$ of ${\Bbb C}[\alpha]$ determined by 
\begin{equation}
s_i(\alpha_i)=-\alpha_i,\ \ 
s_i(\alpha_j)=\alpha_j+\alpha_i\ \ (i\ne j), \ \ \pi(\alpha_j)=\alpha_{j+1}
\end{equation}
for $i,j=0,1,2$. 
When we consider the action of $\widetilde{W}$ on the parameter space $V$, 
we will use the action such that $(w.\varphi)(\bv)=\varphi(w^{-1}.\bv)$ 
for any $\bv\in V$ and $\varphi\in {\Bbb C}[\alpha]$.  
The action of $\widetilde{W}$ on $V$ is given as follows: 
\begin{equation} 
\begin{align}
& s_0.\bv=(v_3+1,v_2,v_1-1),\quad
s_1.\bv=(v_2,v_1,v_3),\\ 
&s_2.\bv=(v_1,v_3,v_2),\quad
\pi.\bv=(v_3+\frac{2}{3},v_1-\frac{1}{3},v_2-\frac{1}{3}).\notag
\end{align}
\end{equation}
for any $\bv=(v_1,v_2,v_3)\in V$.
\par\medskip
One advantage of our representation \eqref{p4f} is that the 
action of the affine Weyl group $\widetilde{W}$ on the fourth Painlev\'e 
equation can be described in a completely symmetric way 
on the dependent variables $f_0$, $f_1$ and $f_2$. 
The action of $\widetilde{W}$ on ${\Bbb C}[\alpha]$ extends in fact to 
the whole differential field $K={\Bbb C}(\alpha;f)$  as follows.

\begin{thm}\label{BT}
The fourth Painlev\'e equation \eqref{p4f} is invariant under the 
following transformations $s_0,s_1,s_2$ and $\pi$: 
\begin{equation}
\begin{aligned}
s_0(f_0)&=f_0, \\
s_0(f_1)&=f_1-{\alpha_0 \over f_0}, \\
s_0(f_2)&=f_2+{\alpha_0 \over f_0}, \\
s_0(\alpha_0)&=-\alpha_0, \\
s_0(\alpha_1)&=\alpha_1+\alpha_0, \\
s_0(\alpha_2)&=\alpha_2+\alpha_0,
\end{aligned}
\ 
\begin{aligned}
s_1(f_1)&=f_1, \\
s_1(f_2)&=f_2-{\alpha_1 \over f_1}, \\
s_1(f_0)&=f_0+{\alpha_1 \over f_1}, \\
s_1(\alpha_1)&=-\alpha_1, \\
s_1(\alpha_2)&=\alpha_2+\alpha_1, \\
s_1(\alpha_0)&=\alpha_0+\alpha_1,
\end{aligned}
\ 
\begin{aligned}
s_2(f_2)&=f_2, \\
s_2(f_0)&=f_0-{\alpha_2 \over f_2}, \\
s_2(f_1)&=f_1+{\alpha_2 \over f_2}, \\
s_2(\alpha_2)&=-\alpha_2, \\
s_2(\alpha_0)&=\alpha_0+\alpha_2, \\
s_2(\alpha_1)&=\alpha_1+\alpha_2,
\end{aligned}
\ \ 
\begin{aligned}
\pi(f_0)&=f_1, \\
\pi(f_1)&=f_2, \phantom{{\alpha_1 \over f}}\\
\pi(f_2)&=f_0,\phantom{{\alpha_1 \over f}} \\
\pi(\alpha_0)&=\alpha_1, \\
\pi(\alpha_1)&=\alpha_2, \\
\pi(\alpha_2)&=\alpha_0.
\end{aligned}
\label{transf}
\end{equation}
Furthermore, these transformations define a representation of 
the affine Weyl group $\widetilde{W}=\langle s_0,s_1,s_2,\pi\rangle$.  
Namely, $\widetilde{W}$ acts on the differential field 
$K={\Bbb C}(\alpha;f)$ as a group of differential automorphisms. 
\end{thm}
\medskip\noindent
Theorem \ref{BT} is proved by straightforward computations.  
The transformations described above will be called the 
{\it B\"acklund transformations} of the fourth Painlev\'e 
equation \eqref{p4f}. 
Note that the independent variable $x=(f_0+f_1+f_2)/3$ is fixed 
under the action of $\widetilde{W}$. 
\par\medskip
Note that, for any $w\in W$, 
one obtains three linear functions 
$\beta_0=w(\alpha_0)$, $\beta_1=w(\alpha_1)$, 
$\beta_2=w(\alpha_2)$ in $\alpha_0,\alpha_1,\alpha_2$. 
Theorem \ref{BT} then implies that, one can specify certain 
rational functions 
$g_0=w(f_0), g_1=w(f_1), g_2=w(f_2)$ in 
$f_0,f_1,f_2,\alpha_0,\alpha_1,\alpha_2$ such that 
\begin{equation}
g_i' + g_i(g_{i+1}-g_{i+2})=\beta_i\quad(i=0,1,2). 
\end{equation} 
Namely, if $(f_0,f_1,f_2)$ is a (generic) solution of \eqref{p4f} with 
parameters $(\alpha_0,\alpha_1,\alpha_2)$, then  
$(g_0, g_1, g_2)$ is again a solution of the same system with 
parameters $(\beta_0,\beta_1,\beta_2)$. 
We give an example below to show how the dependent variables 
$f_0$, $f_1$, $f_2$ are transformed under the action of the 
affine Weyl group.  
\par\medskip
\noindent
{\em Example.}\quad
For $w=s_1s_0$, the B\"acklund transformation 
$w(f_1)=s_1s_0(f_1)$ is computed as follows:
\begin{equation}
f_1 \ {\overset {s_0.} \longrightarrow } \ 
\frac{f_0f_1-\alpha_0}{f_0} \  
{\overset {s_1.} \longrightarrow } \ 
\frac{f_1(f_0f_1-\alpha_0)}
{f_0f_1+\alpha_1}. 
\end{equation}
Similarly we have
\begin{equation}
\begin{align}
&\beta_0=w(\alpha_0)=\alpha_2-3,\quad
\beta_1=w(\alpha_1)=\alpha_0,\quad
\beta_2=w(\alpha_2)=\alpha_1+3; \notag\\
&g_0=w(f_0)=\frac{f_0f_1+\alpha_1}{f_1},\quad
g_1=w(f_1)=\frac{f_1(f_0f_1-\alpha_0)}{f_0f_1+\alpha_1},\\
&g_2=w(f_2)=
\frac{(f_0f_1+\alpha_1)(f_1f_2-\alpha_1)+(3-\alpha_2)f_1^2}
{f_1(f_0 f_1+\alpha_1)}.\notag
\end{align}
\end{equation}
If we specialize these formula to the particular solution 
\begin{equation}\label{ratsol1}
(\alpha_0,\alpha_1,\alpha_2\,;\, f_0,f_1,f_2) =(1,1,1;\,x,x,x), 
\end{equation}
we obtain another rational solution
\begin{equation}\label{ratsol2}
(\alpha_0,\alpha_1,\alpha_2\,;\, f_0,f_1,f_2)= 
(-2,1,4 \,;\, \frac{x^2+1}{x},\frac{x(x^2-1)}{(x^2+1)},
\frac{x^4+2 x^2-1}{x(x^2+1)}).
\end{equation}
A complete description of rational functions in $x$ arising in this way 
will be given later in this paper. 

\par\medskip \noindent
{\em Remark.}\quad
The B\"acklund transformation 
$s_0(f_1)=f_1-\frac{\alpha_0}{f_0}$, for example, becomes singular
when applied to a particular solution such that $f_0=0$. 
This sort of problem, however, is only apparent since such a solution 
arises only under the condition $\alpha_0=0$ as one sees immediately 
from \eqref{p4f}.  
When $\alpha_0=0$, it is natural to understand 
that the B\"acklund transformation $s_0$ 
becomes the identity transformation.   
In general, 
each $g_i=w(f_i)$ is a rational function in $(\alpha; f)$ and 
its denominator possibly becomes identically zero when 
one specializes $(\alpha;f)$ to certain particular solutions.  
Such a phenomenon occurs however only when some of 
the parameters $\alpha_0$, $\alpha_1$, $\alpha_2$ are 
in $3{\Bbb Z}$. 
In such cases, critical factors in the denominator of 
$g_i=w(f_i)$ can actually be eliminated by specializing the parameters 
$(\alpha_0,\alpha_1,\alpha_2)$ in advance.  
With this {\em regularization}, 
our B\"acklund transformations $w(f_i)$ make sense for any particular 
solution. 
\par\medskip
%%%%%%%%%%%%%%%%%%%%%%%%%%%%%%%%%%%%%%%%%%%%%%%%%%%%%%%%%%%%%%%
\section{$\tau$-Functions}
%%%%%%%%%%%%%%%%%%%%%%%%%%%%%%%%%%%%%%%%%%%%%%%%%%%%%%%%%%%%%%%
In this section, we show that our equation \eqref{p4f} for $f=(f_0,f_1,f_2)$ 
can be bilinearized by introducing a triple of $\tau$-functions 
$\tau=(\tau_0,\tau_1,\tau_2)$. 
We also study the B\"acklund transformations on the level of $\tau$-functions.

We introduce the $\tau$-{\em functions} $\tau_0$, $\tau_1$, $\tau_2$ to be 
the dependent variables satisfying the following equations: 
\begin{align}
\label{tau}
f_0&=(\log{\tau_1 \over \tau_2})'+x
=\frac{\tau_1'}{\tau_1}-\frac{\tau_2'}{\tau_2}+x,\notag\\
f_1&=(\log{\tau_2\over \tau_0})'+x
=\frac{\tau_2'}{\tau_2}-\frac{\tau_0'}{\tau_0}+x,\\
f_2&=(\log{\tau_0 \over \tau_1})'+x
=\frac{\tau_0'}{\tau_0}-\frac{\tau_1'}{\tau_1}+x.\notag
\end{align}
We fix the freedom of overall multiplication by a function
in defining $\tau_0,\tau_1,\tau_2$, by imposing the equation
\begin{equation}
2 (\log \tau_0 \tau_1 \tau_2)''+(f_0-x)^2+(f_1-x)^2+(f_2-x)^2=0. 
\end{equation} 
To be more precise, we first introduce a variable $g$ (determined from 
$f_0,f_1,f_2$ up to an additive constant) as an integral of the 
equation
\begin{equation}
\label{intg}
2g'+(f_0-x)^2+(f_1-x)^2+(f_2-x)^2=0.
\end{equation}
Then we require that the $\tau$-functions $\tau_0$, $\tau_1$, $\tau_2$ 
should satisfy 
\begin{equation}
\label{g} 
g=(\log\tau_0\tau_1\tau_2)'=\frac{\tau_0'}{\tau_0}+
\frac{\tau_1'}{\tau_1}+\frac{\tau_1'}{\tau_1}.
\end{equation}
Note that, under the conditions \eqref{tau} and \eqref{g}, 
the $\tau$-functions $\tau_0$,$\tau_1$,$\tau_2$ are  
determined by the equations
\begin{equation}
\begin{align}
\label{taubyf}
(\log\tau_0)'=\frac{\tau_0'}{\tau_0}=\frac{1}{3}(g-f_1+f_2),\notag\\
(\log\tau_1)'=\frac{\tau_1'}{\tau_1}=\frac{1}{3}(g-f_2+f_0),\\
(\log\tau_2)'=\frac{\tau_2'}{\tau_2}=\frac{1}{3}(g-f_0+f_1),\notag
\end{align}
\end{equation}
up to multiplicative constants, respectively.
We remark that the integration constant in $g$ has the effect of 
multiplying each $\tau_i$ by the exponential of a linear function in $x$. 

In order to describe the differential equations to be satisfied by the 
$\tau$-functions, we recall the definition of Hirota's bilinear equations. 
Let $P(\partial_x)$ ($\partial_x=d/dx$) be a linear differential operator
in the $x$-variable with constant coefficients. 
Then Hirota's bilinear operator $P(D_x)$ is defined by 
\begin{equation}
P(D_x) \ F(x)\cdot G(x) = P(\partial_y) F(x+y)G(x-y)\vert_{y=0}, 
\end{equation}
for a given pair of functions $F(x)$, $G(x)$. 

\begin{thm}\label{Hiro}
The fourth Painlev\'e equation \eqref{p4f} for $f_0$, $f_1$, $f_2$,
together with the integral $g$ of \eqref{intg}, is equivalent to the following 
system of Hirota bilinear equations for the triple of $\tau$-functions 
$\tau_0$,$\tau_1$,$\tau_2$:
\begin{align}\label{p4tau}
&(D_x^2-x D_x-\frac{\alpha_0-\alpha_1}{3} ) \ \tau_0 \cdot\tau_1=0, \notag\\ 
&(D_x^2-x D_x-\frac{\alpha_1-\alpha_2}{3} ) \ \tau_1 \cdot\tau_2=0, \\ 
&(D_x^2-x D_x-\frac{\alpha_2-\alpha_0}{3} ) \ \tau_2 \cdot\tau_0=0. \notag
\end{align}
\end{thm}
\medskip\noindent
{\em Proof.} \quad
Note first that, in terms of the logarithms $F_i=\log\tau_i$ ($i=0,1,2$) 
of $\tau$-functions, the dependent variables $f_0,f_1,f_2$ are expressed as follows: 
\begin{equation}
\label{fbyF}
\begin{align}
&f_0=F_1'-F_2'+x, \quad f_1=F_2'-F_0'+x, 
\quad f_2=F_0'-F_1'+x, \\
&g=F_0'+F_1'+F_2'.\notag
\end{align}
\end{equation}
The three equations of Theorem are rewritten into the following equations for $F_0$,
$F_1$, $F_2$:  
\begin{equation}
\label{HinF}
\begin{align}
F_0''+F_1''+(F_0'-F_1')^2-x(F_0'-F_1')-\frac{\alpha_0-\alpha_1}{3}=0,\notag\\
F_1''+F_2''+(F_1'-F_2')^2-x(F_1'-F_2')-\frac{\alpha_1-\alpha_2}{3}=0,\\
F_2''+F_0''+(F_2'-F_0')^2-x(F_2'-F_0')-\frac{\alpha_2-\alpha_0}{3}=0.\notag
\end{align}
\end{equation}
Taking the sum of these three equations, we have 
\begin{equation}
\label{Fsum}
2(F_0''+F_1''+F_2'')+(F_1'-F_2')^2+(F_2'-F_0')^2+(F_0'-F_1')^2=0,
\end{equation}
which corresponds to the equation \eqref{intg} for $g$. 
By subtracting the third equation of \eqref{HinF} from the first, we have
\begin{equation}
F_1''-F_2''-(F_1'-F_2'+x)(2F_0'-F_1'-F_2')-\alpha_0+1=0,
\end{equation}
which corresponds to the differential equation for $f_0$. 
Similarly we have the equations for $f_1$ and $f_2$ from \eqref{HinF}. 
It is also clear that the equations \eqref{HinF} are recovered 
from \eqref{Fsum} and the three equations which correspond to \eqref{p4f}. 
{\qed}
\par\medskip\noindent
{\em Remark.} \quad
Consider the differential field $K(g)={\Bbb C}(\alpha;f)(g)$ obtained 
from $K={\Bbb C}(\alpha;f)$ by adjoining a variable $g$ 
on which the derivation $'$ acts by the formula \eqref{intg}. 
Then Theorem \ref{Hiro} implies that this differential field is 
isomorphic to the differential field ${\Bbb C}(\alpha)(x,F_0',F_1',F_2')$
defined by the relations \eqref{HinF}. 
Note that, by \eqref{HinF} and \eqref{Fsum}, each second 
derivative $F_i''$ 
($i=0,1,2$) can be expressed  in terms of $x$ and $F_0'$, $F_1'$, 
$F_2'$: 
\begin{equation}\label{F2}
F_i''+x(F_{i+1}'-F_{i+2}')+(F_i'-F_{i+1}')(F_i'-F_{i+2}')+
\frac{\alpha_{i+1}-\alpha_{i+2}}{3}=0
\end{equation}
for $i=0,1,2$.
This system is also equivalent to the equation \eqref{p4tau} for 
the triple $\tau_0,\tau_1,\tau_2$ of $\tau$-functions. 
Note that the differential field of our $\tau$-functions is defined as 
${\Bbb C}(\alpha)(x,\tau_0,\tau_1,\tau_2,\tau_0',\tau_1',\tau_2')$ 
by \eqref{F2}, regarded as equations for $\tau$-functions. 
\par\medskip
One important fact is that the action of the affine Weyl group 
on the $f$-variables lifts to the level of $\tau$-functions. 
\begin{thm}\label{BTtau}
The $\tau$-functions $(\tau_0,\tau_1,\tau_2)$ allow 
an action of the affine Weyl group $\widetilde{W}$ which is compatible 
with the action of $\widetilde{W}$ on $f_0,f_1,f_2$ of Theorem \ref{BT}. 
Their B\"acklund transformations are again expressed by Hirota's 
bilinear operators as follows: 
\begin{equation}\label{bktau}
\begin{align}
&s_0(\tau_0)=\frac{1}{\tau_0}(D_x+x)\,\tau_1\cdot\tau_2=
\frac{1}{\tau_0}(\tau_1'\tau_2-\tau_1\tau_2'+x\tau_1\tau_2),\notag\\
&s_1(\tau_1)=\frac{1}{\tau_1}(D_x+x)\,\tau_2\cdot\tau_0=
\frac{1}{\tau_1}(\tau_2'\tau_0-\tau_2\tau_0'+x\tau_2\tau_0),\\
&s_2(\tau_2)=\frac{1}{\tau_2}(D_x+x)\,\tau_0\cdot\tau_1=
\frac{1}{\tau_2}(\tau_0'\tau_1-\tau_0\tau_1'+x\tau_0\tau_1),\notag\\
& s_i(\tau_j)=\tau_j\quad(i\ne j),\quad \pi(\tau_j)=\tau_{j+1}\quad(i,j=0,1,2),\notag
\end{align}
\end{equation}
while $s_0,s_1,s_2$ and $\pi$ act on $\alpha_0$,$\alpha_1$,$\alpha_2$
in the same way as in Theorem \ref{BT}. 
\end{thm}
\medskip\noindent
{\em Proof.} \quad
We first extend the action of $\widetilde{W}$ on ${\Bbb C}(\alpha;f)$ 
to ${\Bbb C}(\alpha;f)(g)$, or equivalently to 
${\Bbb C}(\alpha)(x,F_0',F_1',F'_2)$. 
{}From \eqref{intg} we have 
\begin{equation}
s_0(g')=g' +(f_1-f_2)\frac{\alpha_0}{f_0}-(\frac{\alpha_0}{f_0})^2
=g'-\alpha_0 \frac{f_0'}{f_0^2}
\end{equation}
by \eqref{p4f}. 
Hence we have 
\begin{equation}
s_i(g')=g'-\alpha_i \frac{f_i'}{f_i^2}\quad(i=0,1,2), \quad \pi(g')=g'.
\end{equation}
In view of these, we define the action of $\widetilde{W}$ on $g$ by 
\begin{equation}\label{song}
s_i(g)=g+\frac{\alpha_i}{f_i}\quad (i=0,1,2), \quad \pi(g)=g. 
\end{equation}
One can check that \eqref{song} gives rise in fact to a representation 
of $\widetilde{W}$ on ${\Bbb C}(\alpha;f)(g)$. 
On the variables $F_0', F_1', F_2'$, equation \eqref{taubyf} 
together with \eqref{song} immediately implies 
\begin{equation}
s_i(F_j')=F_j' \quad(i\ne j),\quad \pi(F_j')=F_{j+1}'\quad(i,j=0,1,2). 
\end{equation}  
These formulas justify the definitions of \eqref{bktau} 
other than those for $s_i(\tau_i)$ ($i=0,1,2$). 
As to $s_0(\tau_0)$, we compute 
\begin{equation}
s_0(F_0')=F_0' + \frac{\alpha_0}{f_0}=F_0' + \frac{f_0'}{f_0}+f_1-f_2
=-F_0'+F_1'+F_2'+\frac{f_0'}{f_0}. 
\end{equation}
This leads to the definition
\begin{equation}\label{sandtf}
s_0(\tau_0)=\frac{\tau_1\tau_2}{\tau_0}f_0
=\frac{\tau_1\tau_2}{\tau_0}(\frac{\tau_1'}{\tau_1}-\frac{\tau_2'}{\tau_2}
+x)
=\frac{1}{\tau_0}(D_{x}+x)\,\tau_1\cdot \tau_2. 
\end{equation}
One can check by straightforward computations that 
the definition \eqref{bktau} thus obtained 
defines an action of $\widetilde{W}$ on the differential field
${\Bbb C}(\alpha)(x,\tau_0,\tau_1,\tau_2,\tau_0',\tau_1',\tau_2')$
as a group of differential automorphisms. 
$\qed$
\par\medskip\noindent
We remark that the B\"acklund transformations $s_i(\tau_i)$ 
of Theorem \ref{BTtau} possibly  become zero 
for solutions reducible to Riccati equations, 
while they can be applied repeatedly as long as the $\tau$-functions 
remain nonzero. 
If $(\tau_0,\tau_1,\tau_2)$ is a {\em generic\,} solution, 
we obtain the B\"acklund 
transformations $(w(\tau_0),w(\tau_1),w(\tau_2))$ for 
any $w\in\widetilde{W}$, by Theorem \ref{BTtau}.  

\par
\begin{figure}
\setlength{\unitlength}{1.2mm}
\begin{picture}(50,35)(-33,0)
\hbox{
\hskip4mm
\put(17,0){$s_2(\tau_2)$}
\put(9,15){$\tau_0$}
\put(29,15){$\tau_1$}
\put(-3,30){$s_1(\tau_1)$}
\put(19,30){$\tau_2$}
\put(37,30){$s_0(\tau_0)$}
\put(19,12){$f_2$}
\put(10,22){$f_1$}
\put(26,22){$f_0$}
}
\put(15,15){\vector(1,0){10}}
\multiput(5,30)(20,0){2}{\vector(1,0){10}}
\put(28,13){\vector(-2,-3){6}}
\multiput(18,28)(20,0){2}{\vector(-2,-3){6}}
\put(18, 3){\vector(-2,3){6}}
\multiput( 8,18)(20,0){2}{\vector(-2,3){6}}
\end{picture}
\caption{Six $\tau$-functions}
\label{fig:sixtau}
\end{figure}
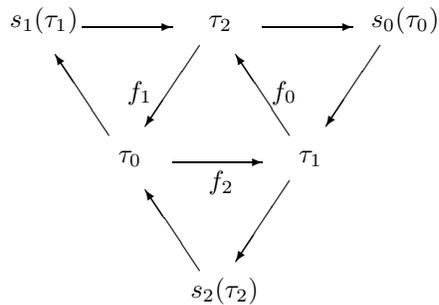

\medskip
{}From the formula \eqref{sandtf} in the proof of Theorem \ref{BTtau}, we have
\begin{cor}\label{fmul}
In terms of the $\tau$-functions $\tau_0$, $\tau_1$, $\tau_2$, 
the dependent variables $f_0$, $f_1$,$f_2$ 
of the fourth Painlev\'e equation \eqref{p4f} are 
expressed multiplicatively as follows: 
\begin{equation}\label{fbyt}
f_0={\tau_0\, s_0(\tau_0) \over \tau_1 \, \tau_2},\quad
f_1={\tau_1\, s_1(\tau_1) \over \tau_2 \, \tau_0},\quad
f_2={\tau_2\, s_2(\tau_2) \over \tau_0 \, \tau_1}. 
\end{equation}
\end{cor} 
\medskip\noindent
The relation between the $f$-variables and the six $\tau$-functions 
in Corollary \ref{fmul} can be represented graphically as in Figure 
\ref{fig:sixtau}. 
Note also that \eqref{fbyt} implies 
\begin{equation}
\tau_0^2 s_0(\tau_0)+\tau_1^2 s_1(\tau_1)+\tau_2^2 s_2(\tau_2)=
3x\, \tau_0 \tau_1 \tau_2 .
\end{equation}
%%%%%%%%%%%%%%%%%

\par\medskip\noindent 
{\em Example.}  As to the rational solution \eqref{ratsol1}
the corresponding $\tau$-functions and their adjacent B\"acklund 
transformations are given by
\begin{equation}
(\tau_0,\tau_1,\tau_2)=(1,1,1), \quad
(s_0(\tau_0),s_1(\tau_1),s_2(\tau_2))=
(x,x,x). 
\end{equation}
As to the rational solution \eqref{ratsol2}, we have 
\begin{equation}
\begin{align}
&(\tau_0,\tau_1,\tau_2)=(x^2+1,x,1), \\
&(s_0(\tau_0),s_1(\tau_1),s_2(\tau_2))=
(1,x^2-1,x^4+2 x^2-1).  \notag
\end{align}
\end{equation} 
These are examples of {\em Okamoto polynomials} which will be discussed 
in the next section. 
\par\medskip
Another corollary of Theorem \ref{BTtau} is the Toda equations 
for our $\tau$-functions. 
\begin{cor}\label{Toda0}
The fourth Painlev\'e equation \eqref{p4tau} for the triple of 
$\tau$-functions $\tau_0,\tau_1,\tau_2$ implies the following 
equation of Toda type: 
\begin{equation}
(\log \tau_0)''+x^2+\frac{\alpha_1-\alpha_2}{3}
=\frac{s_1(\tau_1) s_2(\tau_2)}{\tau_0^2},
\end{equation}
namely,
\begin{equation}
(\frac{1}{2}D_x^2 +x^2+\frac{\alpha_1-\alpha_2}{3})\,\tau_0\cdot\tau_0
=s_1(\tau_1) s_2(\tau_2).
\end{equation}
\end{cor}
\par\noindent
{\em Proof.}  \quad From Corollary \ref{fmul}, we have
\begin{equation}\label{fttt}
f_1 f_2 =\frac{s_1(\tau_1)\,s_2(\tau_2)}{\tau_0^2}. 
\end{equation}
On the other hand, substitution of the formulas
\begin{equation}
F'_2-F_0'=f_1-x,\ \  F'_0-F'_1=f_2-x,\ \ F_1'-F_2'=2x-f_1-f_2 
\end{equation}
into \eqref{F2} with $i=0$ gives 
%\begin{equation}
%F_0''+\frac{x}{3}(\frac{2x}{3}-f_1-f_2)-(f_1-\frac{x}{3})(f_2-\frac{x}{3})
%+\frac{\alpha_1-\alpha_2}{3}=0, 
%\end{equation}
%hence 
\begin{equation}\label{F2f}
F_0''+x^2+\frac{\alpha_1-\alpha_2}{3}=f_1f_2. 
\end{equation}
Equating \eqref{fttt} and \eqref{F2f} we obtain the equation of 
Corollary as desired.
$\qed$. 
\par\medskip\noindent
{\em Remark.} \quad
Our $\tau$-functions are slightly different 
from those introduced by K. Okamoto \cite{Ok}.  
In our formulation, the $\tau$-function in the spirit of Okamoto, 
say $\tau^{\hbox{ok}}$, 
can be defined through the integral of a ``Hamiltonian'' as follows: 
\begin{equation}
H=\frac{1}{3}(f_0 f_1 f_2+\alpha_1 f_2-\alpha_2 f_1)=(\log\tau^{\hbox{ok}})'.
\end{equation}
Note that this implies
$(\log\tau^{\hbox{ok}})''=H'=f_1 f_2.$
Let us introduce the triple of $\tau$-functions of Okamoto type by 
\begin{equation}
(\log \tau_0^{\hbox{ok}})'=H_0, \quad 
(\log \tau_1^{\hbox{ok}})'=H_1, \quad 
(\log \tau_2^{\hbox{ok}})'=H_2, 
\end{equation}
where we define $H_0=H$, $H_1=\pi(H)$, $H_2=\pi^2(H)$ 
by rotation.  
%Then we have 
%\begin{equation}
%(\log \tau_0^{\hbox{ok}})''=f_1 f_2, \quad
%(\log \tau_1^{\hbox{ok}})''=f_2 f_0, \quad
%(\log \tau_2^{\hbox{ok}})''=f_0 f_1. 
%\end{equation}
This implies 
\begin{equation}
f_0=(\log {\tau_1^{\hbox{ok}} \over \tau_2^{\hbox{ok}}})'+\alpha_0 x, \ \ 
f_1=(\log {\tau_2^{\hbox{ok}} \over \tau_0^{\hbox{ok}}})'+\alpha_1 x, \ \ 
f_2=(\log {\tau_0^{\hbox{ok}} \over \tau_1^{\hbox{ok}}})'+\alpha_2 x. 
\end{equation}
(Compare these formulas with our definition \eqref{tau}.)
{}From \eqref{F2f} we also see that 
\begin{equation}
\tau_0^{\hbox{ok}}= e^{x^4/12+(\alpha_1-\alpha_2){x^2}/{6}} \, \tau_0
\end{equation}
up to the multiplication by the exponential of a linear function in $x$. 
\par\medskip
%
%%%%%%%%%%%%%%%%%%%%%%%%%%%%%%%%%%%%%%
\section{Rational solutions}
%%%%%%%%%%%%%%%%%%%%%%%%%%%%%%%%%%%%%%
In this section, we give an explicit description of the rational solutions 
of the fourth Painlev\'e equation \eqref{p4f} in terms of Schur functions.  
\par\medskip
Before discussing the rational solutions, we introduce a family of $\tau$-functions
$(\tau_{m, n})_{m,n\in{\Bbb Z}}$ for the fourth Painlev\'e 
equation \eqref{p4tau}. 
A similar treatment of the lattice of $\tau$-functions has been 
given by K.\,Okamoto \cite{Ok1}. 
We consider the elements 
\begin{equation}
T_1=\pi s_2 s_1, \quad T_2=s_1 \pi s_2
\end{equation}
of the (extended) affine Weyl group $\widetilde{W}$.
Note that these $T_1$ and $T_2$ represent the following 
parallel translations in the parameter space $V$ respectively: 
\begin{equation}
T_1.\bv=\bv+(\frac{2}{3},-\frac{1}{3},-\frac{1}{3}),\quad 
T_2.\bv=\bv+(-\frac{1}{3},\frac{2}{3},-\frac{1}{3}),
\end{equation}
for $\bv\in V$. 
For the triple of $\tau$-functions $(\tau_0,\tau_1, \tau_2)$ of the fourth 
Painlev\'e equation \eqref{p4tau}, 
we introduce an infinite family of dependent variables 
$\tau_{m,n}$ ($m,n\in{\Bbb Z}$) as the B\"acklund transformations 
\begin{equation}
\tau_{m,n} =T_1^m T_2^n (\tau_0)\quad(m,n\in{\Bbb Z}). 
\end{equation}
Note that 
\begin{equation}
T_1(\tau_0)=\tau_1, \quad T_2(\tau_0)=s_1(\tau_1)\quad
\mbox{and} \quad T_2T_1(\tau_0)=T_2(\tau_1)=\tau_2. 
\end{equation}
By these formulas, we have
\begin{equation}
\tau_{0,0}=\tau_0,\quad \tau_{1,0}=\tau_1,\quad
\tau_{1,1}=\tau_2,\quad \tau_{0,1}=s_1(\tau_1). 
\end{equation}
The triple of $\tau$-functions $(\tau_0,\tau_1,\tau_2)$ 
is transformed into 
$(\tau_{m,n},\tau_{m+1,n},\tau_{m+1,n+1})$
by $T_1^m T_2^n$, and into 
$(\tau_{m,n},\tau_{m,n+1},\tau_{m+1,n+1})$
by $T_1^m T_2^n s_1$, respectively. 
The following propositions are obtained 
immediately from the results of the previous section,  
by using the action of $\widetilde{W}$. 
%%%%%%%%%%%%%%%%%%%%%%%%%%%
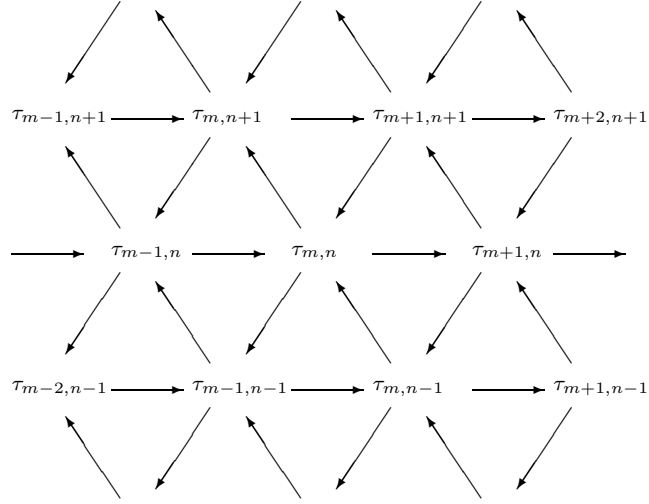
\begin{figure}
\setlength{\unitlength}{1.2mm}
\begin{picture}(80,55)(-25,-10)
\hbox{
\hskip4mm
\put(-4,0){\small$\tau_{m-2,n-1}$}
\put(16,0){\small$\tau_{m-1,n-1}$}
\put(36,0){\small$\tau_{m,n-1}$}
\put(56,0){\small$\tau_{m+1,n-1}$}
\put(7,15){\small$\tau_{m-1,n}$}
\put(27,15){\small$\tau_{m,n}$}
\put(47,15){\small$\tau_{m+1,n}$}
\put(-4,30){\small$\tau_{m-1,n+1}$}
\put(16,30){\small$\tau_{m,n+1}$}
\put(36,30){\small$\tau_{m+1,n+1}$}
\put(56,30){\small$\tau_{m+2,n+1}$}
}
\multiput(7,0)(20,0){3}{\vector(1,0){8}}
\multiput(-4,15)(20,0){4}{\vector(1,0){8}}
\multiput(7,30)(20,0){3}{\vector(1,0){8}}
\multiput(18,-2)(20,0){3}{\vector(-2,-3){6}}
\multiput(8,13)(20,0){3}{\vector(-2,-3){6}}
\multiput(18,28)(20,0){3}{\vector(-2,-3){6}}
\multiput(8,43)(20,0){3}{\vector(-2,-3){6}}
\multiput(8,-12)(20,0){3}{\vector(-2,3){6}}
\multiput(18,3)(20,0){3}{\vector(-2,3){6}}
\multiput( 8,18)(20,0){3}{\vector(-2,3){6}}
\multiput(18,33)(20,0){3}{\vector(-2,3){6}}
\end{picture}
\caption{$\tau$-Functions on the $A_2$-lattice}
\label{fig:taulat}
\end{figure}
%%%%%%%%%%%%%%%%%%%%%%%

\begin{prop}\label{Prop41}
$(1)$ 
For any $m,n\in{\Bbb Z}$, 
the triples 
\begin{equation}\label{twotriples}
(\tau_{m,n},\tau_{m+1,n},\tau_{m+1,n+1})\quad\mbox{and}\quad
(\tau_{m,n},\tau_{m,n+1},\tau_{m+1,n+1})
\end{equation}
 represent the B\"acklund transformations of $(\tau_0,\tau_1,\tau_2)$ for the 
parameters 
\begin{equation}
\begin{align}
&(\alpha_0+3m,\alpha_1+3(n-m),\alpha_2-3n)\quad \mbox{and }\\
&(\alpha_0+\alpha_1+3n,-\alpha_1+3(m-n),\alpha_1+\alpha_2-3m),\notag
\end{align}
\end{equation}
respectively. 
\newline
$(2)$ The corresponding $f$-variables are given respectively by
\begin{equation}
\begin{align}
&\left(\frac{\tau_{m,n}\tau_{m+2,n+1}}{\tau_{m+1,n}\tau_{m+1,n+1}},
\frac{\tau_{m+1,n}\tau_{m,n+1}}{\tau_{m+1,n+1}\tau_{m,n}}, 
\frac{\tau_{m+1,n+1}\tau_{m,n-1}}{\tau_{m,n}\tau_{m+1,n}}
\right)\quad\mbox{and}\\
&\left(\frac{\tau_{m,n}\tau_{m+1,n+2}}{\tau_{m,n+1}\tau_{m+1,n+1}},
\frac{\tau_{m,n+1}\tau_{m+1,n}}{\tau_{m+1,n+1}\tau_{m,n}},
\frac{\tau_{m+1,n+1}\tau_{m-1,n}}{\tau_{m,n}\tau_{m,n+1}}
\right).\notag
\end{align}
\end{equation}
\end{prop}
\begin{prop}\label{Prop41b}
$(1)$ The family of $\tau$-functions $\tau_{m,n}$ $(m,n\in {\Bbb Z})$  
satisfies the following three types of bilinear equations: 
\begin{equation}
\begin{align}
&(D_x+x)\, 
\tau_{m,n}\cdot\tau_{m+1,n}=\tau_{m,n-1}\tau_{m+1,n+1},\notag\\
&(D_x+x)\, 
\tau_{m,n}\cdot\tau_{m,n+1}=\tau_{m+1,n+1}\tau_{m-1,n},\\
&(D_x+x)\, 
\tau_{m,n}\cdot\tau_{m-1,n-1}=\tau_{m-1,n}\tau_{m,n-1}.\notag
\end{align}
\end{equation}
$(2)$ The family of $\tau$-functions $\tau_{m,n}$ $(m,n\in {\Bbb Z})$  
satisfies the following three types of Toda equations: 
\begin{equation}
\begin{align}
&(\frac{1}{2}D_x^2+x^2-\frac{2\alpha_1+\alpha_2}{3}+2m-n)\, 
\tau_{m,n}\cdot\tau_{m,n}=\tau_{m+1,n}\tau_{m-1,n},\notag\\
&(\frac{1}{2}D_x^2+x^2+\frac{\alpha_1-\alpha_2}{3}-m+2n)\, 
\tau_{m,n}\cdot\tau_{m,n}=\tau_{m,n+1}\tau_{m,n-1},\\
&(\frac{1}{2}D_x^2+x^2+\frac{\alpha_1+2\alpha_2}{3}-m-n)\, 
\tau_{m,n}\cdot\tau_{m,n}=\tau_{m-1,n-1}\tau_{m+1,n+1}.\notag
\end{align}
\end{equation}
\end{prop}
\medskip\noindent
%We remark that, 
%in the coordinates $(v_1,v_2,v_3)$ as in \eqref{parameter-alpha}, 
%the triples of \eqref{twotriples} correspond to the triples of $\tau$-functions 
%for the parameters
%\begin{equation}
%\begin{align}
%&(v_1+\frac{-2m+n}{3},v_2+\frac{m-2n}{3},v_3+\frac{m+n}{3})\quad\mbox{and}\\
%&(v_2+\frac{m-2n}{3},v_1+\frac{-2m+n}{3},v_3+\frac{m+n}{3})\notag
%\end{align}
%\end{equation}
%respectively.
%A graphical description of the Painlev\'e system on 
%the $A_2$-lattice will be given in Section 6. 
%\par\medskip\noindent
{\em Remark.}\quad
As we already remarked in the previous section, B\"acklund transformations for 
$\tau$-functions possibly become singular, when applied to particular solutions 
which are reducible to Riccati equations. 
In such cases, 
we need to restrict the indices $(m,n)$ for $\tau_{m,n}$ 
to a region of ${\Bbb Z}^2$ bounded by certain lines on which $\tau_{m,n}=0$. 
\par\medskip 

%%%%%%%%%%%%%%%%%%%%%%%%%%%%%%%%%%%%%%%%%%%%
All the rational solutions of \eqref{p4f} are obtained from 
\begin{equation}
\begin{align} \label{ratsolAB}
\mbox{(A)} &\qquad(\alpha_0,\alpha_1,\alpha_2; f_0,f_1,f_2) 
= (1,1,1;\, x,x,x), \quad\mbox{or} \\
\mbox{(B)} &\qquad(\alpha_0,\alpha_1,\alpha_2; f_0,f_1,f_2) = (3,0,0; 3x,0,0).\notag
\end{align}
\end{equation}
by B\"acklund transformations.  
We will determine the $\tau$-functions $\tau_{m,n}$($m,n\in{\Bbb Z}$) 
for these rational solutions.  
\par\medskip
In the case of B\"acklund transformations of the rational 
solution (A) of \eqref{ratsolAB}, 
the $\tau$-functions $\tau_{m,n}$ ($m,n\in{\Bbb Z}$) turn out to be polynomials, 
which we call the {\em Okamoto polynomials}.
We recall that the $\tau$-functions for (A) are given by
\begin{equation}\label{A0}
(\alpha_0,\alpha_1,\alpha_2; \tau_0,\tau_1,\tau_2)=
(1,1,1; 1,1,1). 
\end{equation}
\begin{thm} \label{Thm42}
The $\tau$-functions $\tau_{m,n}(x)$ for the solution \eqref{A0}
are polynomials in $x$.  
These polynomials $\tau_{m,n}(x)=Q_{m,n}(x)$ $(m,n\in{\Bbb Z})$
are characterized by the Toda equations
\begin{equation}
\begin{align}
&(\frac{1}{2}D_x^2+x^2-1+2m-n)\, 
Q_{m,n}\cdot Q_{m,n}=Q_{m+1,n}Q_{m-1,n},\notag\\
&(\frac{1}{2}D_x^2+x^2-m+2n)\, 
Q_{m,n}\cdot Q_{m,n}=Q_{m,n+1} Q_{m,n-1},\\
&(\frac{1}{2}D_x^2+x^2+1-m-n)\, 
Q_{m,n}\cdot Q_{m,n}=Q_{m-1,n-1} Q_{m+1,n+1}\notag
\end{align}
\end{equation}
with initial condition
\begin{equation}
Q_{0,0}=Q_{1,0}=Q_{1,1}=1, \quad
Q_{2,1}=x. 
\end{equation}
\end{thm}
\medskip\noindent
We remark that $Q_m(x)=Q_{m,0}(x)$ and $R_{m}(x)=Q_{m+1,1}(x)$ 
($m\in {\Bbb Z}$) are the 
original Okamoto polynomials discussed in \cite{FOU}. 
In fact, they are determined by the recurrence relations
\begin{equation}
(\frac{1}{2}D_x^2+x^2+2m-1)\, 
Q_{m}\cdot Q_{m}=Q_{m+1}Q_{m-1}\quad (m\in{\Bbb Z})
\end{equation}
with initial condition 
$ Q_0=Q_1=1 $, and by 
\begin{equation}
(\frac{1}{2}D_x^2+x^2+2m)\, 
R_{m}\cdot R_{m}=R_{m+1}R_{m-1}\quad(m\in{\Bbb Z})
\end{equation}
with $R_0=1, R_1=x$, respectively.
%%%%%%%%%%%%%%%%%%%%%%%%%%%%%%%%%%%%%%%%%%%%
%\par\medskip
The fact that $\tau_{m,n}(x)$ are polynomials will be proved in 
Section 5 in the course of the proof of Theorem \ref{ThmOka} below. 
The other statements in Theorem \ref{Thm42} are consequences of 
Proposition \ref{Prop41b}.
\par\smallskip
The $\tau$-functions for the rational solution (B) of \eqref{ratsolAB} 
are given by
\begin{equation}\label{B0}
(\alpha_0,\alpha_1,\alpha_2; \tau_0,\tau_1,\tau_2)=
(3,0,0; e^{-x^4/12}, e^{-x^4/12+x^2/2},e^{-x^4/12-x^2/2}). 
\end{equation}
\begin{thm} \label{Thm43}
The $\tau$-functions $\tau_{m,n}(x)$ for the solution \eqref{B0}
are defined for $(m,n)\in{\Bbb Z}^2$ with  $m\ge n\ge0$. 
They can be written in the form
\begin{equation}\label{expHer}
\tau_{m,n}(x)=\exp\left(-\frac{x^4}{12}+\frac{m-2n}{2}x^2\right) H_{m-n,n}
\quad(m\ge n\ge 0), 
\end{equation}  
for some polynomials $H_{m,n}(x)$. 
These polynomials $H_{m,n}(x)$ $(m,n\ge 0)$ are characterized by 
the Toda equations
\begin{equation}
\begin{align}
&(\frac{1}{2}D_x^2+3m)\, 
H_{m,n}\cdot H_{m,n}=H_{m+1,n} H_{m-1,n},\\
&(\frac{1}{2}D_x^2-3n)\, 
H_{m,n}\cdot H_{m,n}=H_{m,n+1}H_{m,n-1},\notag
\end{align}
\end{equation}
with initial condition
\begin{equation}
H_{0,0}=H_{1,0}=H_{0,1}=1\quad \mbox{and}\quad
H_{1,1}=3x.
\end{equation}
\end{thm}
\medskip\noindent
We remark that $H_{m,1}(x)$ and $H_{1,m}$ ($m=0,1,2,\ldots$) coincide with the 
Hermite polynomials up to rescaling. 
We will call $H_{m,n}(x)$ $(m,n\ge 0)$ the {\em generalized Hermite polynomials}. 
The fact that $\tau_{m,n}(x)$ are expressed as in \eqref{expHer} will be 
proved in Section 5 in the course of the proof of Theorem \ref{ThmHer} below.  
\par\medskip
The Okamoto polynomials $Q_{m,n}(x)$ ($m,n\in{\Bbb Z}$) and 
the generalized Hermite polynomials $H_{m,n}(x)$  ($m,n\ge 0 $) 
are in fact expressible in terms of Schur functions. 
We recall the definition of Schur functions in order to make this statement precise. 
\par\medskip
A {\em partition} $\lambda=(\lambda_1,\lambda_2,\ldots)$ 
(or a {\em Young diagram}) is a sequence of non-negative integers such that 
$\lambda_1 \geq \lambda_2 \geq  \cdots \geq 0$
and that $\lambda_i=0$ for $i\gg 0$.
The number of nonzero parts $\lambda_i$ is called the {\em length} of $\lambda$
and denoted by $l(\lambda)$.
For each partition $\lambda$, we define the {\em Schur function} 
$S_{\lambda}(t)=S_{\lambda}(t_1,t_2,\ldots)$ by 
\begin{equation}
S_{\lambda}(t)=\det \big( p_{\lambda_i-i+j}(t) \big)_{
1 \leq i,j \leq l(\lambda)}, 
\end{equation}
where $p_n(t)$ are the polynomials in $t$ determined by the generating 
function
\begin{equation}
\exp \left(\sum_{k=1}^{\infty} t_k z^k \right)=
\sum_{n=0}^{\infty}p_n(t)z^n. 
\end{equation}
(We set $p_n(t)=0$ for $n<0$.)  
Note that  $p_n(t)$ can be defined equivalently by 
\begin{equation}
p_n(t)={\sum}_{k_1+2k_2+\cdots+n k_n=n} \ 
\frac{t_1^{k_1} t_2^{k_2}\cdots t_n^{k_n}}{k_1! k_2!\cdots k_n !}. 
\end{equation}

We say that a subset $M\subset{\Bbb Z}$ is a {\em Maya diagram} if 
\begin{equation}
m \in M \quad (m \ll 0)
\quad\mbox{and}\quad
m \not\in M \quad (m \gg 0).
\end{equation}
To each Maya diagram 
$M=\{\ldots, m_3,m_2,m_1\} \quad
(\cdots < m_3 < m_2 < m_1)$, 
one can associate a unique partition $\lambda=(\lambda_1,\lambda_2,\ldots)$ 
such that $m_i-m_{i+1}=\lambda_i-\lambda_{i+1}+1$ for $i=1,2,\ldots$. 
Note that all the Maya diagrams  
$M+k=\{\ldots,m_2+k,m_1+k\}$ ($k\in{\Bbb Z}$)
obtained from $M=\{\ldots, m_3,m_2,m_1\}$ by shifting 
define the same partition by this correspondence. 
For each pair $(m,n)$ of integers, we define the Maya diagram 
$M(m,n)$ as follows: 
\begin{equation}
M(m,n)=3 D_m \cup(3 D_n+1) \cup (3 D_0+2), 
\end{equation}
where
\begin{equation}
D_l=\{ n\in {\Bbb Z} \ | \ n<l \}\quad(l\in{\Bbb Z}). 
\end{equation}
We denote by  $\lambda(m,n)$ 
the partition corresponding to $M(m,n)$. 
Partitions of the form $\lambda(m,n)$ ($m,n\in{\Bbb Z})$ are 
called the 3-{\em reduced partitions}. 
We remark that a partition $\lambda$ is 3-reduced if and only 
if $\lambda$ has no hook with length of a multiple of $3$.  
Also, Schur functions $S_{\lambda(m,n)}(t)$ for 3-reduced partitions 
are called {\em 3-reduced Schur functions}.  
It is known that a Schur function $S_{\lambda}(t)$ is 3-reduced 
if and only if 
\begin{equation}
\partial_{t_{3n}} S_{\lambda}(t)=0 \quad \mbox{for all }\quad n=1,2,\cdots.
\end{equation}
\par\medskip
\begin{thm} \label{ThmOka}
Each Okamoto polynomial $Q_{m,n}(x)$ $(m,n\in{\Bbb Z})$ is a monic 
polynomial of degree $m^2+n^2-mn-m$ with integer coefficients.   
It is expressed by the 3-reduced Schur function as 
\begin{equation}
Q_{m,n}(x)=N_{m,n} S_{\lambda(m,n)}(x,\frac{1}{2},0,0,\ldots), 
\end{equation}
where $N_{m,n}$ is a positive integer determined by the hook-length formula. 
\end{thm}
%\FigOka
%\ExOka 
\par\noindent
{\em Example.}\quad 
The Maya diagram $M(3,2)$ is obtained from $D_3,D_2,D_0$ as follows.
\par\medskip\smallskip\noindent
\begin{picture}(250,105)(-60,-5)
\hbox{
\multiput(62.5,82.5)(10,0){8}{$\bullet$} \put(24,82){\small $D_3$}
\multiput(62.5,72.5)(10,0){7}{$\bullet$} \put(24,72){\small $D_2$}
\multiput(62.5,62.5)(10,0){5}{$\bullet$} \put(24,62){\small $D_0$}
\multiput(44,63)(0,10){3}{$\cdots$}
\put(98,93){\small $-1$}
\put(113,93){\small $0$}
\put(123,93){\small $1$}
\put(133,93){\small $2$}
\put(143,93){\small $\cdots$}
\put(83,93){\small $\cdots$}
\put(0,8){$\implies$}
\put(90,32){\small $M(3,2)$}
\put(34,7.5){\small $\cdots$}
\multiput(48.5,7.5)(10,0){6}{$\bullet$}
\put(118.5,7.5){$\bullet$}
\put(128.5,7.5){$\bullet$}
\put(148.5,7.5){$\bullet$}
\put(43,17){\small $-3$}
\put(74,17){\small $-1$}
\put(89,17){\small $0$}
\put(99,17){\small $1$}
\put(109,17){\small $2$}
\put(119,17){\small $3$}
\put(129,17){\small $4$}
\put(139,17){\small $5$}
\put(149,17){\small $6$}
\put(159,17){\small $\cdots$}
}
\multiput(40,60)(0,10){4}{\line(1,0){150}}
\multiput(60,60)(10,0){12}{\line(0,1){30}}
\put(110,55){\line(0,1){45}}
\multiput(30,5)(0,10){2}{\line(1,0){170}}
\multiput(46,5)(10,0){15}{\line(0,1){10}}
\multiput(56,0)(30,0){5}{\line(0,1){25}}
\end{picture}
\comment{
\par\medskip\smallskip\noindent
\begin{picture}(250,35)(-15,-5)
\hbox{
\multiput(41.5,17.5)(8,0){8}{$\bullet$} \put(4,18){\small $D_3$}
\multiput(41.5,9.5)(8,0){7}{$\bullet$} \put(4,10){\small $D_2$}
\multiput(41.5,1.5)(8,0){5}{$\bullet$} \put(4,2){\small $D_0$}
\multiput(24,1.5)(0,8){3}{$\cdots$}
\put(68,26){\small $-1$}
\put(82,26){\small $0$}
\put(90,26){\small $1$}
\put(98,26){\small $2$}
\put(106,26){\small $\cdots$}
\put(52,26){\small $\cdots$}
\put(152,10){\small $\implies$}
\put(180,24){\small $M(3,2)$}
\put(184,5.5){\small $\cdots$}
\multiput(197.5,5.5)(8,0){6}{$\bullet$}
\put(253.5,5.5){$\bullet$}
\put(261.5,5.5){$\bullet$}
\put(277.5,5.5){$\bullet$}
\put(192,14){\small $-3$}
\put(217,14){\small $-1$}
\put(230,14){\small $0$}
\put(238,14){\small $1$}
\put(246,14){\small $2$}
\put(254,14){\small $3$}
\put(262,14){\small $4$}
\put(270,14){\small $5$}
\put(278,14){\small $6$}
\put(286,14){\small $\cdots$}
}
\multiput(20,0)(0,8){4}{\line(1,0){120}}
\multiput(40,0)(8,0){12}{\line(0,1){24}}
\put(80,-4){\line(0,1){36}}
\multiput(180,4)(0,8){2}{\line(1,0){140}}
\multiput(196,4)(8,0){15}{\line(0,1){8}}
\multiput(204,0)(24,0){5}{\line(0,1){20}}
\end{picture}
} % \endcomment
\par\medskip\noindent
Hence we have $M(3,2)=\{\ldots,-2,-1,0,1,3,4,6\}$ and 
\begin{equation}
\lambda(3,2)=(2,1,1)={\Young(2,1,1)}.
\end{equation}
In this case, the Schur function $S_{(2,1,1)}(t)$ and the Okamoto polynomials 
$Q_{3,2}(x)$ are 
\begin{equation}
\begin{align} 
S_{(2,1,1)}(t)&=\frac{1}{8}t_1^4- \frac{1}{2} t_1^2 t_2 - \frac{1}{2} t_2^2 + t_4,
\\
Q_{3,2}(x)&=x^4-2x^2-1,\notag
\end{align}
\end{equation}
respectively. 
A typical sequence of 3-reduced partitions is given by 
\begin{equation}
\lambda(m,0)=(2m-2,2m-4,\ldots,2)\quad \mbox{for}\quad m>0,
\end{equation}
which corresponds to the Okamoto polynomials $Q_m(x)$  for $m>0$. 
For other examples, see Figure \ref{fig:okamoto} in the next section. 

\par\medskip
The generalized Hermite polynomials $H_{m,n}(x)$ 
are expressed by the Schur functions for rectangular Young diagrams 
$\lambda=(n^m)=(n,n,\ldots,n,0,0,\ldots)$. 
\begin{thm}\label{ThmHer}
Each generalized Hermite polynomial $H_{m,n}(x)$ $(m,n\ge 0)$ is 
a polynomial of degree $mn$ with rational coefficients. 
It can be written as 
\begin{equation}
H_{m,n}(x)=C_{m,n} S_{(n^m)}(x,\frac{1}{6},0,0,\ldots),
\end{equation}
where the normalization constant is given by
\begin{equation} 
C_{m,n}=(-1)^{n(n-1)/2} 3^{(m+n)(m+n-1)/2}(m+n-1)^!
\end{equation}
with $n^!=n!(n-1)!(n-2)!\cdots 2! 1!$.
\end{thm}
\noindent
The relationship between the sequences $H_{1,n}(x)$, $H_{m,1}(x)$ 
and the ordinary Hermite polynomials $H_n(x)$ is obvious since 
\begin{equation}
\begin{align}
&\sum_{n=0}^\infty \, S_{(n)}(x,\frac{1}{6},0,0,\ldots) z^n =
\exp\left(x z +\frac{1}{6} z^2 \right),\\
&\sum_{m=0}^\infty \, S_{(1^m)}(x,\frac{1}{6},0,0,\ldots) z^m =
\exp\left(x z -\frac{1}{6} z^2 \right),\notag
\end{align}
\end{equation}
while the Hermite polynomials have the generating function
\begin{equation}
\sum_{n=0}^\infty \, \frac{1}{n!} H_n(x) z^n = \exp\left(2xz - z^2 \right).
\end{equation}
%\FigHer
%\ExHer
\par\medskip
The proof of Theorems \ref{ThmOka} and \ref{ThmHer} will be given in the 
next section. 

%%%%%%%%%%%%%%%%%%%%%%%%%%%%%%%%%%%%%%
\section{Proof of Theorems  \ref{ThmOka} and \ref{ThmHer}}
%%%%%%%%%%%%%%%%%%%%%%%%%%%%%%%%%%%%%%
The Hirota bilinear equations for our $\tau$-functions (Theorem \ref{Hiro})
arise naturally from the so-called {\em modified KP hierarchy}
\cite{JM} by certain similarity reduction (see also \cite{NT}, \cite{QNC}). 
This fact is the key to the proof of Theorems \ref{ThmOka} and \ref{ThmHer}. 
\par\medskip
Consider two functions $G_0(t_1,t_2)$ and $G_1(t_1,t_2)$ in 
the two variables $(t_1,t_2)$, and suppose that they 
satisfy the following Hirota bilinear equation
\begin{equation}
\label{mKP0}
(D_{t_1}^2+D_{t_2}) \, G_0(t_1,t_2) \cdot G_1(t_1,t_2) =0. 
\end{equation}
Introducing the degrees of $t_1,t_2$ by $\deg t_1=1$ and 
$\deg t_2=2$,
we assume that each $G_i$ is homogeneous of degree $d_i\in {\Bbb C}$ 
for $i=0,1$:
\begin{equation}\label{sim}
(t_1\partial_{t_1}+2 t_2\partial_{t_2}) G_i(t_1,t_2)=d_i G_i(t_1,t_2). 
\end{equation}
Fixing a constant $k\in{\Bbb C}^\times$, 
we define the functions $\tau_0(x),\tau_1(x)$ in one variable by
\begin{equation}
\tau_i(x)=G_i(x,k)\quad (k\in {\Bbb C^\times}; \ i=0,1). 
\end{equation}
Formally,  each $G_i$ is recovered by the formula
\begin{equation}
G_i(t_1,t_2)= \left(\frac{t_2}{k}\right)^{d_i \over 2} 
\tau_i\left(\frac{t_1}{\sqrt{t_2/k}}\right).
\end{equation}
Then it is easy to check 
\begin{lem}\label{s-reduction}\label{Lem51}
Under the similarity condition \eqref{sim}, the equation \eqref{mKP0}
for the pair $G_0(t_1,t_2)$, $G_1(t_1,t_2)$ is equivalent to the Hirota 
bilinear equation 
\begin{equation}
\left(2k D_x^2 -xD_x+d_0-d_1\right)\ \tau_0(x)\cdot\tau_1(x)=0, 
\end{equation}
for $\tau_0(x)$, $\tau_1(x)$. 
\end{lem}
{}From this lemma with $k=1/2$, we immediately have

\begin{prop}\label{simred}
The fourth Painlev\'e equation \eqref{p4tau} for 
the triple of $\tau$-functions 
$\tau_0(x), \tau_1(x), \tau_2(x)$ 
is equivalent to the similarity reduction of the Hirota equations
\begin{equation}
(D_{t_1}^2+D_{t_2} ) \ G_{i}(t_1,t_2) \cdot G_{i+1}(t_1,t_2)=0 
\quad (i=0,1,2)
\end{equation}
for three functions $G_i(t_1,t_2)$ $(i=0,1,2)$ in two variables. 
The similarity condition is given by 
\begin{equation}
G_{i}(t_1,t_2)=(2t_2)^{d_i/2} \, 
\tau_i\left(\frac{t_1}{\sqrt{2t_2}} \right)
\quad(i=0,1,2),
\end{equation}
and the parameters are related by 
\begin{equation}\label{abyd}
\begin{align}
&\alpha_0=1-2d_0+d_1+d_2, \notag \\
&\alpha_1=1+d_0-2d_1+d_2, \\
&\alpha_2=1+d_0+d_1-2d_2. \notag
\end{align}
\end{equation}
\end{prop}

\par\medskip
Recall that the (first) {\em modified  KP hierarchy} \cite{JM}
is the following system of Hirota bilinear equations 
for a pair of $\tau$-functions 
$\tau_0(t)$  and $\tau_1(t)$ in infinite time 
variables $t=(t_1,t_2,\ldots)$: 
\begin{equation}\label{mKP}
\sum_{n=0}^\infty p_n(-2s)\,p_{n+2}(\widetilde{D}_t)
\exp\left(\sum_{m=1}^\infty s_m D_{t_m} \right)\ 
\tau_{0}(t) \cdot \tau_{1}(t)=0, 
\end{equation}
where $s=(s_1,s_2,\ldots)$ are parameters 
and $\widetilde{D}_t=(D_{t_1}/1,D_{t_2}/2,\ldots)$. 
The constant term of \eqref{mKP} with respect to $s$ 
implies the bilinear equation
\begin{equation}
(D_{t_1}^2+D_{t_2})\  \tau_0(t)\cdot \tau_{1}(t)=0, 
%&(D_{t_1}^3-4 D_{t_3}-3D_{t_1}D_{t_2})\ 
%\tau_\ell(t)\cdot \tau_{\ell+1}(t)=0,\notag\\
%&(D_{t_1}^4+8 D_{t_1}D_{t_3}+3D_{t_1}^2D_{t_2}-6D_{t_2}^2)\ 
%\tau_\ell(t)\cdot \tau_{\ell+1}(t)=0,\notag\\
%& \cdots\notag
\end{equation}
which is nothing but the equation \eqref{mKP0} discussed above. 
For the proof of Theorems \ref{ThmOka} and \ref{ThmHer}, 
we will recall the following fact about Schur functions 
from the theory of KP hierarchy. 
Let $X_m=X_m(t;\partial_t)$ ($m\in {\Bbb Z}$) be the 
{\em vertex operators} of the KP hierarchy defined by the generating 
function
\begin{equation}
X(z)=\sum_{m\in {\Bbb Z}} X_m z^m=
\exp\left(\sum_{k=1}^\infty t_k z^k\right)
\exp\left(-\sum_{k=1}^\infty \frac{z^{-k}}{k}\partial_{t_k} \right). 
\end{equation}
Then we have 
\begin{lem}\label{Lem52}
For any partition $\lambda$ and $k \in {\Bbb Z}$, the pair
$\tau_0(t)=S_{\lambda}(t)$ and  
$\tau_1(t)=X_k S_{\lambda}(t)$ solves the first modified KP
hierarchy \eqref{mKP}.
In particular we have
\begin{equation}
(D_{t_1}^2+D_{t_2})\  \tau_0(t)\cdot \tau_1(t)=0.
\end{equation}
\end{lem}
\medskip\noindent
We will give a proof of this lemma in Appendix for completeness. 
\par\medskip
All the Schur functions $S_{\lambda}(t)$ are obtained 
from $S_\emptyset(t)=1$ by applying vertex operators
repeatedly:
\begin{equation}\label{SinX}
S_{\lambda}(t)=
X_{\lambda_1}\cdots X_{\lambda_n}.1,
\end{equation}
for any partition $\lambda=(\lambda_1,\ldots,\lambda_n,0,\ldots)$.
The action of vertex operators on Schur functions can be 
computed by \eqref{SinX} together with the commutation relations
\begin{equation}\label{comXX}
X_kX_l=-X_{l-1}X_{k+1}, \quad X_k.1=0\ \ (k<0),\quad X_0.1=1,
\end{equation}
where $k, l\in{\Bbb Z}$. (See Appendix.) 
A more systematic way is to use Maya diagrams. 
For a given Maya diagram $M$, 
let $\lambda$ be the corresponding partition and suppose that
$l(\lambda)\le n$. 
Then we have 
\begin{equation}\label{XS}
X_k.S_\lambda(t)=
\left\{
\begin{array}{cc}
\pm S_\mu(t)\quad&\mbox{if}\quad k+n\notin M,\\
0\quad&\mbox{if}\quad k+n\in M,
\end{array}\right.
\end{equation}
for each $k\in{\Bbb Z}$. 
Here $\mu$ stands for  the partition corresponding to the 
Maya diagram $M\cup\{k+n\}$.  The sign in this formula 
is determined by the parity of the number of integers 
$m\in M$ such that $m>k+n$. 
%%%%%%%%%%%%%%%%%%%%%%%%%%%%%%%%%%%%%%%
\par\medskip
\noindent
{\em Proof of Theorem \ref{ThmOka} for Okamoto polynomials.}\quad 
By using the formula \eqref{XS}, one can compute how the 
3-reduced Schur functions are transformed by vertex operators. 
\begin{lem}\label{Lem53}
As to the action of the vertex operators, 
we have the following two types of cyclic relations among 
3-reduced Schur functions
\begin{align}\label{cyc1}
X_{2m-n}\,.S_{\lambda(m,n)}(t)&=\pm S_{\lambda(m+1,n)}(t),\notag\\
X_{2n-m}\,.S_{\lambda(m+1,n)}(t)&=\pm S_{\lambda(m+1,n+1)}(t),\\
X_{-m-n}\,.S_{\lambda(m+1,n+1)}(t)&=\pm S_{\lambda(m,n)}(t),\notag
\end{align}
and 
\begin{align}
X_{2n-m+1}\,.S_{\lambda(m,n)}(t)&=\pm S_{\lambda(m,n+1)}(t),\notag\\
X_{2m-n-1}\,.S_{\lambda(m,n+1)}(t)&=\pm S_{\lambda(m+1,n+1)}(t),\\
X_{-m-n}\,.S_{\lambda(m+1,n+1)}(t)&=\pm S_{\lambda(m,n)}(t),\notag
\end{align}
for any $m,n\in{\Bbb Z}$. 
\end{lem}
\par\noindent
{\em Example.}\quad
Consider the 3-reduced partitions 
$\lambda(3,1)=(3,1)$, $\lambda(4,1)=(5,3,1)$ and 
$\lambda(4,2)=(4,2,1,1)$. 
%By \eqref{SinX} and \eqref{comXX}, we compute
%\begin{equation}
%\begin{align}
%X_{-1}.S_{(5,3,1)}&=X_{-1}X_{5}X_{3}X_{1}.1=-X_{4}X_{0}X_{3}X_{1}.1 \\
%&=X_{4}X_{2}X_{1}X_{1}.1=S_{(4,2,1,1)}.\notag
%\end{align}
%\end{equation}
%Hence we have $X_{-1}.S_{(5,3,1)}=S_{(4,2,1,1)}$.  
For this triple of partitions, 
we have a `cycle' of 3-reduced Schur functions 
%%%%%%%%%%%%%%%%%
\par\noindent
{
\setlength{\unitlength}{1.2mm}
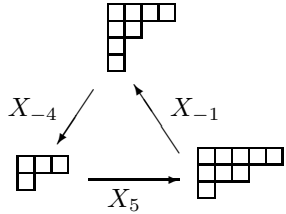
\begin{picture}(50,27)(-33,10)
\hbox{
\hskip4mm
\put(9,15){$\Young(3,1)$}
\put(29,15){$\Young(5,3,1)$}
\put(19,30){$\Young(4,2,1,1)$}
\put(19,12){$X_5$}
\put(8,22){$X_{-4}$}
\put(26,22){$X_{-1}$}
}
\put(17,15){\vector(1,0){10}}
\put(17.5,25){\vector(-2,-3){4}}
\put( 27,18){\vector(-2,3){5}}
\end{picture}
}
\par\noindent
which is an example of \eqref{cyc1} for $(m,n)=(3,1)$. 
Notice also that the index of each vertex operator $X_k$ represents the 
difference of degrees of Schur functions. 
\par\medskip
{}From Lemma \ref{Lem53} together with Lemma \ref{Lem52}, 
we obtain two types of triples of 
3-reduced Schur functions satisfying the bilinear equations 
of Proposition \ref{simred}. 
Note that the 3-reduced Schur function $S_{\lambda(m,n)}$ is 
homogeneous of degree 
\begin{equation}
d_{m,n}=|\lambda(m,n)|=m^2+n^2-mn-m
\end{equation}
with respect to the degree defined by $\deg t_i=i$ ($i=1,2,\ldots$).  
Set
\begin{equation}
s_{m,n}(x)=S_{\lambda(m,n)}(x,\frac{1}{2},0,0,\ldots)
\end{equation}
for any $m,n\in {\Bbb Z}$. 
Then, by combining Lemma \ref{Lem53} and 
Proposition \ref{simred}, we have

\begin{prop}\label{PropOka}
$(1)$  For any $m,n\in{\Bbb Z}$, the triple 
\begin{equation}
(s_{m,n}(x), s_{m+1,n}(x), s_{m+1,n+1}(x))
\end{equation}
solves the fourth Painlev\'e equation \eqref{p4tau} with the 
parameters 
\begin{equation}
(\alpha_0,\alpha_1,\alpha_2)=(3m+1,3(n-m)+1,-3 n+1).
\end{equation}
$(2)$  For any $m,n\in{\Bbb Z}$, the triple 
\begin{equation}
(s_{m,n}(x), s_{m,n+1}(x), s_{m+1,n+1}(x))
\end{equation}
solves the fourth Painlev\'e equation \eqref{p4tau} with 
parameters 
\begin{equation}
(\alpha_0,\alpha_1,\alpha_2)=(3n+2,3(m-n)-1,-3 m+2).
\end{equation}
\end{prop}
\medskip\noindent
We remark that, in the coordinates $(v_1,v_2,v_3)$ of the 
parameter space $V$ as in \eqref{parameter-alpha}, the triples of 
$\tau$-functions in this proposition give rise to 
solutions with parameters 
\begin{equation}
%(v_1,v_2,v_3)=(\frac{1-2m+n}{3},\frac{m-2n}{3},\frac{-1+m+n}{3})
(v_1,v_2,v_3)=(\frac{1}{3},0,-\frac{1}{3})
-m(\frac{2}{3},-\frac{1}{3},-\frac{1}{3})
-n(-\frac{1}{3},\frac{2}{3},-\frac{1}{3})
\end{equation}
and
\begin{equation}
%(v_1,v_2,v_3)=(\frac{m-2n}{3},\frac{1-2m+n}{3},\frac{-1+m+n}{3}),
(v_1,v_2,v_3)=(0,\frac{1}{3},-\frac{1}{3})
-m(-\frac{1}{3},\frac{2}{3},-\frac{1}{3})
-n(\frac{2}{3},-\frac{1}{3},-\frac{1}{3}),
\end{equation}
respectively. 

It is clear that each triple of $\tau$-functions of Proposition 
\ref{PropOka}
defines a rational solution of \eqref{p4f} in $f$-variables.  
For each $(\alpha_0,\alpha_1,\alpha_2)$ of this proposition, 
the fourth Painlev\'e equation \eqref{p4f} has 
a unique rational solution by \cite{M}.  
Hence we conclude that each $s_{m,n}(x)$ is a constant 
multiple of the $\tau$-function $\tau_{m,n}(x)$ 
for the solution (A) of \eqref{ratsolAB}. 
This shows that $\tau_{m,n}(x)$ are in fact polynomials 
in $x$.   
The assertion that the Okamoto polynomials 
$\tau_{m,n}(x)=Q_{m,n}(x)$ are monic polynomials with 
integer coefficients follows 
either from the B\"acklund transformations
or from the Toda equations of Proposition \ref{Prop41b}. 

%%%%%%%%%%%%%%%%%%%%%%%%%%%%%%%%%%%%%%%%%%%%%%
\par\medskip
\noindent
{\em Proof of Theorem \ref{ThmHer} for generalized Hermite 
polynomials.} 
\quad 
By the formula \eqref{XS}, we have 
\begin{lem}
Under the action of vertex operators, we have the following 
relations among Schur functions for rectangular Young diagrams: 
\begin{equation}
\begin{align}
X_{-m} . S_{((n+1)^m)}(t)&=(-1)^m S_{(n^m)}(t), \notag \\
X_{n-m} . S_{((n+1)^m)}(t)&=(-1)^m S_{(n^{(m+1)})}(t), \\ 
X_n . S_{(n^m)}(t)&=S_{(n^{(m+1)})}(t), \notag 
\end{align}
\end{equation}
for $m,n\ge 0$. 
\end{lem}

For each $m,n\ge 0$, let 
\begin{equation}
h_{m,n}(x)=S_{(n^m)}(x,{1 \over 6},0,0,\ldots)
\end{equation}
be the specialization of the Schur function 
$S_{\lambda}(t)=S_{(n^m)}(t)$ 
associated with rectangular Young diagram $\lambda=(n^m)$.
Then by Lemma \ref{s-reduction} (with $k=1/6$), we have
\begin{lem}
\begin{align}
&\left( D_x^2-3 x D_x+3m \right) h_{m,n+1} \cdot h_{m,n}=0, \notag \\
&\left( D_x^2-3 x D_x+3(m-n) \right) h_{m,n+1} \cdot h_{m+1,n}=0, \\
&\left( D_x^2-3 x D_x-3n \right) h_{m,n} \cdot h_{m+1,n}=0. \notag
\end{align}
\end{lem}
These relations do not fit directly for the triple of 
$\tau$-functions as in \eqref{p4tau} 
since they do not make a `cycle'. 
This problem can be repaired however by changing the normalization of 
$h_{m,n}$ as follows: 
\begin{equation}
u_{m,n}(x)=\exp \left(-{x^4 \over 12}+{m-n \over 2}x^2 \right) h_{m,n}(x).
\end{equation}
Then we have 
\begin{prop}\label{PropHer}
\begin{align}\label{ueq}
&\left( D_x^2-x D_x+m+2n+1 \right) u_{m,n+1} \cdot u_{m,n}=0, \notag \\
&\left( D_x^2-x D_x+m-n \right) u_{m+1,n} \cdot u_{m,n+1}=0, \\
&\left( D_x^2-x D_x-2m-n-1 \right) u_{m,n} \cdot u_{m+1,n}=0. \notag
\end{align}
Namely, the triple 
\begin{equation}
(\tau_0,\tau_1,\tau_2)=(u_{m,n},u_{m+1,n},u_{m,n+1})
\end{equation}
solves the fourth Painlev\'e equation \eqref{p4tau} with parameters
\begin{equation}
(\alpha_0,\alpha_1,\alpha_2)=(3(m+n+1),-3 m,-3n).
\end{equation}
\end{prop}
\par\noindent
{\em Proof.} \quad
Note that the Hirota bilinear equations have the following 
formulas of Leibniz type: 
\begin{align}
&D_x(g_1 u_1 \cdot g_2 u_2)=
D_x (g_1 \cdot g_2) u_1 u_2+g_1 g_2 D_x(u_1 \cdot u_2), \notag \\
&D_x^2 (g_1 u_1 \cdot g_2 u_2)=
D_x^2 (g_1 \cdot g_2) u_1 u_2+2 D_x(g_1 \cdot g_2) D_x(u_1 \cdot u_2)+
g_1 g_2 D_x^2(u_1 \cdot u_2).\notag
\end{align}
Applying these to $g_i=\exp(x^4/12+a_i x^2/2)$ $(i=1,2)$, we have 
\begin{align}
&(D_x^2-3x D_x+\beta) (g_1 u_1 \cdot g_2 u_2) \notag \\
=g_1 g_2 &
\big\{D_x^2+(2a_{12}-3)x D_x+
(2-3a_{12}+a_{12}^2)x^2+(a_1+a_2)+\beta
\big\} u_1 \cdot u_2, \notag
\end{align}
where $a_{12}=a_1-a_2$. 
The equations \eqref{ueq} can be checked easily by using this formula. 
$\qed$
\par\medskip\noindent
We remark that, in the coordinates $(v_1,v_2,v_3)$ of $V$, 
the solutions of Proposition \ref{PropHer} have the parameters
\begin{equation}
%(v_1,v_2,v_3)=(\frac{-2m-n}{3},\frac{m-n}{3},\frac{m+2n}{3}).
(v_1,v_2,v_3)=-m(\frac{2}{3},-\frac{1}{3},-\frac{1}{3})+
n(-\frac{1}{3},-\frac{1}{3},\frac{2}{3})\quad(m,n= 0,1,2,\ldots). 
\end{equation}
\par\medskip
As in the case of Okamoto polynomials, 
we see that each $\tau_{m,n}$ ($m\ge n \ge 0$) for the solution (B) 
of \eqref{ratsolAB} is a constant multiple of $u_{m-n,n}$ by 
comparing the parameters.  
Hence we see that $\tau_{m,n}(x)$ has the expression 
of \eqref{expHer}. 
The only problem remaining is to fix the constant factors. 
The leading coefficient of the polynomial
$H_{m,n}(x)$ is given by
\begin{equation}
(-1)^{n(n-1)/2} (m-1)^{!} (n-1)^{!} 3^{(m+n)(m+n-1)/2},
\end{equation}
which can be determined inductively by the Toda equations of 
Theorem \ref{Thm43}. 
On the other hand, the leading coefficient of $h_{m,n}(x)$ is determined 
as  
\begin{equation}
(m-1)^! (n-1)^! /(m+n-1)^!, 
\end{equation}
by the hook-length formula. 
Hence we have Theorem \ref{ThmHer}. 
\par\medskip
%%%%%%  \endcomment 
\par\medskip
We show in Figures \ref{fig:okamoto} and \ref{fig:hermite} below 
how the $\tau$-functions for rational solutions are arranged 
on the $A_2$-lattice.  
Also, we include some examples of Okamoto polynomials and 
generalized Hermite polynomials of small degrees. 
%%%%%%%%%%%%%%%%%%%%%%%%%%%%%%%%%%%%%%%%%%%%
%\noindent
%{\em Example.}
%\begin{align}
%M(3,1)&=\{\ldots,-2,-1,0,1,2,4,7\}, \quad
%\lambda(3,1)=\Young(3,1), \notag \\
%M(4,1)&=\{\ldots,-2,-1,0,1,2,4,7,10\}, \quad
%\lambda(4,1)=\Young(5,3,1), \notag \\
%M(4,2)&=\{\ldots,-2,-1,0,1,2,4,5,7,10\}, \quad
%\lambda(4,2)=\Young(4,2,1,1), \notag
%\end{align}
%
%The Figure \ref{fig:okamoto} shows the way to generate
%all 3-reduced Schur functions by applying Vertex operators.
%For instance, the triple of 3-reduced Schur functions in
%the above Example are related as follows
%\begin{align}
%&X_5 . {\Young(3,1)}={\Young(5,3,1)}, \notag \\
%&X_{-1} . {\Young(5,3,1)}={\Young(4,2,1,1)}, \notag \\
%&X_{-4} . {\Young(4,2,1,1)}={\Young(3,1)}. \notag
%\end{align}
%In general, the number $k$ of $X_k$ is given by the difference
%of degree of the Young diagrams.
%To derive the equations above, we have used the following relations
%\begin{align}
%&S_{(\lambda_1,\ldots,\lambda_n)}(t)=
%X_{\lambda_1}\cdots X_{\lambda_n} . 1, \\
%&X_k X_l=-X_{l-1} X_{k+1}, \quad X_j . 1=0 \ (j <0), \quad 
%X_0 .1=1, \notag
%\end{align}
%where $\lambda$ is a partition and $k,l \in {\Bbb Z}$.
%
%
%%%%%%%%%%%%%%%%%%%%%%%%%%%%%%%%%%%%%%%%%%
\par\newpage
\FigOka
\ExOka
\par\newpage
\FigHer
\ExHer
\par\newpage
%%%%%%%%%%%%%%%%%%%%%%%%%%%%%%%%%%%%%%%%%%
\def\thesection{A}
%%%%%%%%%%%%%%%%%%%%%%%%%%%%%%%%%%%%%%%%
\section{Appendix}
%%%%%%%%%%%%%%%%%%%%%%%%%%%%%%%%%%%%%%%%
%
In this Appendix, we give a brief summary of relevant facts 
on Schur functions and their
relation to KP-hierarchy for the sake of reference. 
\subsection{Schur functions.}\quad
A partition $\lambda=(\lambda_1,\lambda_2,\ldots)$ is a 
sequence of non-negative integers such that 
$\lambda_1 \geq \lambda_2 \geq \cdots \geq 0$
and that $\lambda_i=0$ for $i\gg 0$. 
The number of nonzero $\lambda_i$ is called the length of $\lambda$
and denoted by $l(\lambda)$.
For each partition $\lambda$, the Schur function $S_\lambda(t)=
S_{\lambda}(t_1,t_2,\ldots)$ is defined as follows: 
\begin{equation}
\label{Jac}
S_{\lambda}(t)=
\det \big( p_{\lambda_i-i+j}(t) \big)_{1 \leq i,j \leq l(\lambda)},
\end{equation}
where $p_n(t)$ are the polynomials defined by the generating 
function
\begin{equation}
\exp \left(\sum_{k=1}^{\infty} t_k z^k \right)=
\sum_{n=0}^{\infty}p_n(t)z^n. 
\end{equation}
Usually, the Schur functions are defined as the following 
character polynomials of the general linear group $GL(N,{\Bbb C})$ 
($N \geq l(\lambda)$): 
\begin{equation}
s_{\lambda}(x_1,\ldots,x_N)={
\det(x_j^{\lambda_i+\delta_i}) \over
\det(x_j^{\delta_i})},
\end{equation}
where $\delta_i=N-i$ ($i=1,\ldots,N$).
The polynomials 
$S_{\lambda}(t)$ and $s_{\lambda}(x)$ are related by 
$S_{\lambda}(t)=s_{\lambda}(x)$, where $t_k=\sum_{i=1}^n (x_i^k)/k$.
In this context, the formula \eqref{Jac} above is the Jacobi-Trudi 
formula representing $s_\lambda(x)$ in terms of complete 
homogeneous symmetric functions. 

The coefficients of $S_{\lambda}(t)$ with respect to the $t$-variables
are related with irreducible character $\pi_{\lambda}$ of 
the symmetric group ${\frak S}_n$ of degree 
$n=\vert \lambda \vert=\sum_{i} \lambda_i$ as follows:
\begin{equation}
S_{\lambda}(t)=\sum_{m_1,m_2,\ldots \geq 0}
\pi_{\lambda}(1^{m_1}2^{m_2}\cdots)
{t_1^{m_1} \over m_1 !}
{t_2^{m_2} \over m_2 !} \cdots,
\end{equation}
where $\pi_{\lambda}(1^{m_1}2^{m_2}\cdots)$ is the character value
on the conjugate class of cycle type $(1^{m_1}2^{m_2}\cdots)$.
In particular, the coefficient of $t_1^n$ is given by the 
hook-length formula
\begin{equation}
\frac{\pi_\lambda(1^n)}{n!} = \prod_{s\in\lambda} \frac{1}{h(s)},
\end{equation}
where $h(s)=\lambda_i+\lambda_j'-i-j+1$, $\lambda'$ being 
the conjugate partition, denotes the hook-length of $\lambda$
at $s=(i,j)$.
\subsection{(Modified) KP hierarchy}
In the following, we use the notation
\begin{equation}
\xi(z,t)=\sum_{n=1}^{\infty} t_n z^n, \quad
\xi(z^{-1},\widetilde{\partial}_t)=\sum_{n=1}^{\infty} {z^{-n} \over n}
\partial_{t_n}.
\end{equation}
Consider the operators 
$V_k=V_k(z,t)$ ($k \in {\Bbb Z}$) defined by 
\begin{equation}
V_k=e^{k \xi(z,t)} e^{-k \xi(z^{-1},\widetilde{\partial}_t)}.
\end{equation}
For each $m\in {\Bbb Z}$, 
we define the operators $X_m$ and $X^{\ast}_m$ 
as the coefficient of $z^m$ in $V_1$ and $V_{-1}$,
respectively: 
\begin{align}
V_1(z,t)=&X(z,t)=\sum_{m \in {\Bbb Z}} X_m z^m,\\
V_{-1}(z,t)=&X^{\ast}(z,t)=\sum_{m \in {\Bbb Z}} X^{\ast}_{m} z^m.
\notag
\end{align}
By using the formula
\begin{equation}\label{normal-order}
V_k(z,t) V_l(w,t)= \left( 1-{w \over z} \right)^{k l} 
e^{k \xi(z,t)+l \xi(w,t)}
e^{-k \xi(z^{-1},\widetilde{\partial}_t)-l \xi(w^{-1},\widetilde{\partial}_t)},
\end{equation}
we obtain
\begin{lem}
The vertex operators $X_m$ and $X_m^\ast$ $(m\in{\Bbb Z})$ 
satisfy the following anti-commutation relations: 
\begin{align}
\label{comm}
X_{m} X_{n}+X_{n-1} X_{m+1}&=0, \notag\\
X^{\ast}_{m} X^{\ast}_{n}+
X^{\ast}_{n-1} X^{\ast}_{m+1}&=0, \\
X_{m} X^{\ast}_{n}+
X^{\ast}_{n+1} X_{m-1}&=\delta_{m+n,0}. \notag
\end{align}
\end{lem}
\begin{prop}\label{genS}
For any partition $\lambda=(\lambda_1,\lambda_2,\ldots)$ of 
length $l(\lambda)\le l$, we have
\begin{equation}
S_{\lambda}(t)=X_{\lambda_1} \cdots X_{\lambda_l} .1.
\end{equation}
\end{prop}
 
\noindent{\em Proof.}\quad
By using \eqref{normal-order}, we have
\begin{equation}
X(z_1,t) \cdots X(z_l,t) . 1=
\prod_{1 \leq i < j \leq l} \left(1-{z_j \over z_i} \right)
\prod_{i=1}^{l} \exp \left( \sum_{n=1}^{\infty} t_n z_i^n \right).
\end{equation}
By taking the coefficient of 
$z^{\lambda}=z_1^{\lambda_1}\cdots z_l^{\lambda_l}$
of this expression, we obtain the formula \eqref{Jac}.
$\qed$
\par\medskip
The KP hierarchy is a system of nonlinear partial differential 
equations for an unknown function $\tau(t)=\tau(t_1,t_2,\ldots)$ 
including the Hirota bilinear equation
\begin{equation}
(D_{t_1}^4-4 D_{t_1} D_{t_3}+3 D_{t_2}^2) \tau(t) \cdot \tau(t)=0.
\end{equation}
The whole system of the KP hierarchy is represented by
the following bilinear relation: 
\begin{equation}\label{KPXX}
\oint {dz \over 2 \pi i} 
X^{\ast}(z,t') \tau(t') X(z,t) \tau(t)=0.
\end{equation}
\begin{prop}
For any partition $\lambda$, the 
Schur function $S_\lambda(t)$ 
is a solution of the KP hierarchy.
\end{prop}

\noindent{\em Proof.} \quad 
Note first that the bilinear equation
\eqref{KPXX} can be rewritten as follows: 
\begin{equation}\label{XX}
\left(\sum_{m+n=-1} X^{\ast}_{m} \otimes X_{n} \right)
\tau \otimes \tau=0.
\end{equation}
Here $\tau \otimes \tau=\tau(t') \tau(t)$ is regarded as an element of 
${\Bbb C}[[t']] \otimes {\Bbb C}[[t]]$.
By the anti-commutation relation \eqref{comm}, one has
\begin{align}
&\left( \sum_{m+n=-1} X^{\ast}_{m} \otimes X_{n} \right) \ 
X_{k} \otimes X_{k} \\
&=X_{k+1} \otimes X_{k-1} \
\left( \sum_{m+n=-1} X^{\ast}_{m} \otimes X_{n} \right)-
1 \otimes X_{k-1} X_{k}, \notag
\end{align}
and the last term $X_{k-1} X_{k}$ vanishes.
Hence, by applying the operator $X_{k+1} \otimes X_{k-1}$ to 
\eqref{XX}, it follows that $X_{k} \tau(t)$ is also a solution of 
the KP hierarchy. Starting from the solution
$\tau(t)=1$, we see that 
all the Schur functions are solutions of KP hierarchy
by Proposition \ref{genS}
$\qed$
\begin{prop}
Let $\tau_0(t)=\tau(t)$ be any solution of the KP hierarchy, and put
\begin{equation}
\tau_1(t)=X(w,t) \tau(t).
\end{equation}
Then, we have 
\begin{equation}\label{mKPXX}
\oint {dz \over 2 \pi i} \,z\,
X^{\ast}(z,t') \tau_0(t') \ 
X(z,t) \tau_1(t)=0.
\end{equation}
%Especially,
%\begin{equation}
%(D_1^2+D_2) \tau \cdot \tau_1=0.
%\label{mKP1}
%\end{equation}
\end{prop}
\medskip\noindent
{\em Proof.}\quad
Apply $X(w,t)$ to the second factor 
$X(z,t) \tau(t)$ of the bilinear relation \eqref{KPXX}.
Then one obtains \eqref{mKPXX} by using 
the relation
\begin{equation}
X(w,t) X(z,t) \tau(t)= -{z \over w} X(z,t) X(w,t) \tau(t)=
-{z \over w} X(z,t) \tau_1(t) \notag
\end{equation}
as desired. $\qed$\par
\medskip\noindent
The formula \eqref{mKPXX} is the bilinear relation 
of the first modified KP hierarchy.  
\par\medskip
By the change of variables $t \rightarrow t-s$ and
$t' \rightarrow t+s$, the relations \eqref{KPXX} and \eqref{mKPXX}
can be rewritten into the following systems of 
Hirota bilinear equations 
\begin{align}
&\sum_{n=0}^{\infty} p_{n}(-2 s) p_{n+1}({\widetilde D_t}) 
\exp \left( \sum_{m=1}^{\infty} s_m D_{t_m} \right) \tau(t) \cdot \tau(t)=0, \\
&\sum_{n=0}^{\infty} p_{n}(-2 s) p_{n+2}({\widetilde D_t}) 
\exp \left( \sum_{m=1}^{\infty} s_m D_{t_m} \right) \tau_0(t) \cdot \tau_1(t)=0,
\label{mKPh}
\end{align}
where ${\widetilde D}_{t_n}=D_{t_n} /n$.
These are the Hirota bilinear equations for the $\tau$-functions 
of the KP hierarchy and the first modified KP hierarchy, respectively. 

%%%%%%%%%%%%%%%%%%%%%%%%%%%%%%%%%%%%%%%%

\end{document}